\documentclass[a4paper, 11pt]{amsproc}

\usepackage{fullpage}
\usepackage{amsmath, amssymb, amsfonts, amsthm, amssymb, mathtools, mathrsfs}
\usepackage{ascmac}
\usepackage{comment}
\usepackage{bm}

\usepackage{hyperref}

\usepackage{graphicx}
\usepackage{here}
\usepackage{time}
\usepackage[abbrev]{amsrefs}

\usepackage{xcolor}
\usepackage[capitalize,nameinlink,noabbrev,nosort]{cleveref}
\hypersetup{
	colorlinks=true,       
	linkcolor=blue,          
	citecolor=blue,        
	filecolor=blue,      
	urlcolor=blue,           
}

\makeatletter
\@namedef{subjclassname@2020}{%
  \textup{2020} Mathematics Subject Classification}
\makeatother

\newtheorem{theoremcounter}{Theorem Counter}[section]

\theoremstyle{definition}
\newtheorem{defi}[theoremcounter]{Definition}
\newtheorem{ex}[theoremcounter]{Example}

\theoremstyle{plain}
\newtheorem{lem}[theoremcounter]{Lemma}
\newtheorem{prop}[theoremcounter]{Proposition}
\newtheorem{cor}[theoremcounter]{Corollary}

\newtheorem{thm}[theoremcounter]{Theorem}

\numberwithin{equation}{section}\numberwithin{figure}{section}

\theoremstyle{remark}
\newtheorem{rem}[theoremcounter]{Remark}

\newcommand{\z}{\mathfrak{z}}
\newcommand{\bbH}{\mathbb{H}}
\newcommand{\calI}{\mathcal{I}}
\newcommand{\calL}{\mathcal{L}}

\DeclareMathOperator{\ImNew}{Im}
\renewcommand{\Im}{\ImNew}
\DeclareMathOperator{\ReNew}{Re}
\renewcommand{\Re}{\ReNew}

\DeclareMathOperator{\SL}{SL}
\DeclareMathOperator{\GL}{GL}
\DeclareMathOperator{\sgn}{sgn}
\DeclareMathOperator{\tr}{tr}
\DeclareMathOperator{\val}{val}
\DeclareMathOperator{\arcsinh}{arcsinh}
\DeclareMathOperator{\arccosh}{arccosh}

\newcommand{\smat}[1]{\left(\begin{smallmatrix}#1\end{smallmatrix}\right)}
\newcommand{\pmat}[1]{\begin{pmatrix}#1\end{pmatrix}}

\newcommand{\hyp}[1]{d_{\mathrm{hyp}}(#1)}

\theoremstyle{thmstyleone}%
\theoremstyle{thmstyletwo}%

\theoremstyle{thmstylethree}%

\begin{document}

\title{A Hyperbolic Analogue of the Rademacher Symbol}

\author{Toshiki Matsusaka}
\address{Faculty of Mathematics, Kyushu University, Motooka 744, Nishi-ku, Fukuoka 819-0395, Japan}
\email{matsusaka@math.kyushu-u.ac.jp}


\maketitle

\begin{abstract}
One of the most famous results of Dedekind is the transformation law of $\log \Delta(z)$. After a half-century, Rademacher modified Dedekind's result and introduced an $\SL_2(\mathbb{Z})$-conjugacy class invariant (integer-valued) function $\Psi(\gamma)$ called the Rademacher symbol. Inspired by Ghys' work on modular knots, Duke--Imamo\={g}lu--T\'{o}th (2017) constructed a hyperbolic analogue of the symbol.

In this article, we study their hyperbolic analogue of the Rademacher symbol $\Psi_\gamma(\sigma)$ and provide its two types of explicit formulas by comparing it with the classical Rademacher symbol. In association with it, we contrastively show Kronecker limit type formulas of the parabolic, elliptic, and hyperbolic Eisenstein series. These limits give harmonic, polar harmonic, and locally harmonic Maass forms of weight 2.
\end{abstract}

\setcounter{tocdepth}{2}
\tableofcontents

\section{Introduction} \label{s1}

The \emph{discriminant function}
\[
	\Delta(z) := q \prod_{n=1}^\infty (1 - q^n)^{24}
\]
for $q = e^{2\pi i z}$, $\Im(z) > 0$ is one of the most classical functions in the theory of modular forms. Let $\Gamma := \SL_2(\mathbb{Z})$ be the full modular group, and $\bbH := \{z \in \mathbb{C} \mid \Im(z) > 0\}$ the upper half-plane. It is well-known that the function $\Delta(z)$ satisfies the equation $\Delta (\gamma z) = (cz + d)^{12} \Delta(z)$
for any $\gamma = \smat{a & b \\ c & d} \in \Gamma$ with the M\"{o}bius transformation $\gamma z := (az+b)/(cz+d)$. This equation implies that its holomorphic logarithm $\log \Delta(z)$ also satisfies
\begin{align}\label{R-definition}
	\log \Delta \left(\frac{az+b}{cz+d} \right) - \log \Delta(z) = 12 \sgn(c)^2 \log \left(\frac{cz+d}{i \sgn(c)} \right) + 2\pi i \Phi(\gamma)
\end{align}
for some integer-valued function $\Phi(\gamma)$. Here we fix the meaning of the holomorphic logarithm as
\[
	\log \Delta(z) = 2\pi i z - 24 \sum_{n=1}^\infty \sum_{m=1}^\infty \frac{q^{mn}}{m}
\]
and take the principal branch $\Im \log z \in (-\pi, \pi]$. The symbol $\sgn(c) \in \{-1, 0, 1\}$ is the usual sign function. For $c = 0$, we interpret the first term on the right-hand side as 0. In his paper~\cite{Ded92}, Dedekind studied the function $\Phi: \Gamma \to \mathbb{Z}$, and gave its explicit formula in terms of the Dedekind sum $s(a,c)$. 

We first briefly review Rademacher's works~\cite{Rad56,RG72} on this topic. The \emph{Dedekind sum} is defined by $s(a,c) := \sum_{k=1}^{c-1} ((k/c)) ((ka/c))$ with $((x)) := x - \lfloor x \rfloor - 1/2$ if $x \in \mathbb{R} \setminus \mathbb{Z}$ and zero otherwise. Dedekind established the formula
\[
	\Phi \left( \pmat{a & b \\ c & d} \right) = \begin{cases}
		\dfrac{a+d}{c} - 12 \sgn(c) \cdot s(a, \lvert c \rvert) &\text{if } c \neq 0,\\
		\dfrac{b}{d}  &\text{if } c=0.
	\end{cases}
\]
After a half-century, Rademacher slightly modified the definition of $\Phi$ to get a class invariant function $\Psi (\gamma)$. To be exact, he added one term as
\begin{align}\label{Ded-Rade-rel}
	\Psi(\gamma) := \Phi(\gamma) - 3\sgn(c(a+d)).
\end{align}
Then the function $\Psi(\gamma)$ is a class invariant. More strongly, 
\[
	\Psi(\gamma) = \Psi(-\gamma) = -\Psi(\gamma^{-1}) = \Psi(g^{-1} \gamma g)
\]
holds for any $g \in \Gamma$. We now call $\Psi(\gamma)$ the \emph{Rademacher symbol}. 
Looking at the definitions, the Dedekind sum and the Rademacher symbol seem like specific objects. However, they are very ubiquitous. 
For example, see Asai~\cite{Asa70}, Hirzebruch~\cite{Hir71}, Atiyah~\cite{Ati87}, Arakawa~\cite{Ara88},  Barge--Ghys~\cite{BG92}, Kirby--Melvin~\cite{KM94}, Ghys~\cite{Ghy07}, Zagier~\cite{Zag10}, Li~\cite{Li18}, Kaneko--Mizuno~\cite{KM20}, and so on.

All non-scalar elements $\pm I \neq \gamma \in \Gamma$ are classified into three cases as follows.
\begin{enumerate}
	\item $\gamma$ is \emph{parabolic} if $\lvert \mathrm{tr}(\gamma)\rvert = 2$. Then $\gamma$ has the unique fixed point in $\mathbb{Q} \cup \{i\infty \}$.
	\item $\gamma$ is \emph{elliptic} if $\lvert \mathrm{tr}(\gamma)\rvert < 2$. Then $\gamma$ has a fixed point in $\bbH$.
	\item $\gamma$ is \emph{hyperbolic} if $\lvert \mathrm{tr}(\gamma)\rvert > 2$. Then $\gamma$ has two fixed points in $\mathbb{R}$.
\end{enumerate}
Here the fixed point $z$ of $\gamma$ means the solution in $\bbH \cup \mathbb{R} \cup \{i\infty\}$ of the equation $\gamma z = z$. For hyperbolic elements $\gamma \in \Gamma$, the Rademacher symbol $\Psi(\gamma)$ has some extra properties. We then focus on the case when $\gamma = \smat{a & b \\ c & d} \in \Gamma$ is hyperbolic and exhibit some known results. For simplicity, we assume that $\sgn(c) > 0$ and $\mathrm{tr}(\gamma) > 2$ by replacing $\gamma$ with one of $\pm \gamma^{\pm 1}$ if needed. First, the Rademacher symbol $\Psi(\gamma)$ has an expression in terms of a cycle integral of the weight $2$ mock modular form $E_2(z)$. In fact, by the definition of $\log \Delta(z)$,
\[
	\frac{1}{2\pi i} \frac{d}{dz} \log \Delta(z) = 1 - 24 \sum_{m=1}^\infty \left(\sum_{d\mid m} d \right) q^m =: E_2(z).
\]
A hyperbolic element $\gamma$ has two fixed points $w_\gamma > w'_\gamma$ in $\mathbb{R}$. Let $S_\gamma$ be the geodesic in $\bbH$ connecting $w_\gamma$ and $w'_\gamma$, and take a point $z_0 \in S_\gamma$. Then the integral on the geodesic segment in $S_\gamma$ gives
\[
	\int_{z_0}^{\gamma z_0} E_2(z) dz = \frac{1}{2\pi i} \bigg(\log \Delta(\gamma z_0) - \log \Delta(z_0) \bigg) = \frac{6}{\pi i} \log \left(\frac{cz_0+d}{i} \right) + \Phi(\gamma).
\]
Since $\gamma^n z_0 \to w_\gamma$ as $n \to \infty$ and $cw_\gamma + d > 1$, we get the integral expression of the Rademacher symbol
\begin{align}\label{Eint-Psi}
	\lim_{n \to \infty} \Re \int_{\gamma^n z_0}^{\gamma^{n+1} z_0} E_2(z) dz = \Phi(\gamma) - 3 = \Psi(\gamma).
\end{align}

Second, the Rademacher symbol has an explicit formula in terms of the continued fraction coefficients of $w_\gamma$. Assume that the hyperbolic element $\gamma$ is conjugate to a $\gamma'$ of the form
\[
	\gamma \sim \gamma' = \pmat{a_0 & 1 \\ 1 & 0} \pmat{a_1 & 1 \\ 1 & 0} \cdots \pmat{a_{2n-1} & 1 \\ 1 & 0}
\]
with $a_j \geq 1$. The $\gamma'$ has the fixed point $w_\gamma = [\overline{a_0, a_1, \dots, a_{2n-1}}]$ whose continued fraction expansion is purely periodic (see \cref{s2}). Then Lang~\cite[(2.17)]{Lan76} and Zagier~\cite[Section V, Lemma]{Zag75} established the following simple formula (see also Kaneko--Mizuno~\cite[Appendix 2]{KM20}).
\begin{align}\label{Psi-explicit}
	\Psi(\gamma) = \sum_{j=0}^{2n-1} (-1)^j a_j.
\end{align}

In 2007, Ghys~\cite{Ghy07} gave a topological interpretation to the Rademacher symbol $\Psi(\gamma)$ for a hyperbolic $\gamma$. It is well-known that the quotient $\SL_2(\mathbb{Z}) \backslash \SL_2(\mathbb{R})$ is diffeomorphic to the knot complement of a trefoil knot $K$ in the $3$-sphere $S^3$. For each hyperbolic element $\gamma \in \Gamma$, he introduced a knot $C_\gamma$ in $\Gamma \backslash \SL_2(\mathbb{R}) \cong S^3 \backslash K$ as follows. Suppose that a hyperbolic element $\gamma = \smat{a & b \\ c & d} \in \Gamma$ satisfies the above assumptions $c > 0, a+d > 2$, and is primitive, that is, there is no $\gamma' \in \Gamma$ and $n > 1$ such that $\gamma = \pm \gamma'^n$. The $\gamma$ is diagonalized by the scaling matrix
\begin{align}\label{scaling-mat}
	M_\gamma:= \frac{1}{\sqrt{w_\gamma - w'_\gamma}} \pmat{w_\gamma & w'_\gamma \\ 1 & 1} \in \SL_2(\mathbb{R}),
\end{align}
that is,
\begin{align}\label{fund-unit}
	M_\gamma^{-1} \gamma M_\gamma = \pmat{c w_\gamma + d & 0 \\ 0 & cw'_\gamma + d} =: \pmat{\xi_\gamma & 0 \\ 0 & \xi_\gamma^{-1}}
\end{align}
with $\xi_\gamma > 1$. Then a curve
\[
	C_\gamma(t) := M_\gamma \pmat{e^t & 0 \\ 0 & e^{-t}}, \qquad 0 \leq t \leq \log \xi_\gamma
\]
on the geodesic flow defines a simple closed curve in $\Gamma \backslash \SL_2(\mathbb{R})$. We call it a \emph{modular knot}. It is then a natural question to ask whether we can compute a knot invariant of the modular knot $C_\gamma$ in terms of $\gamma$. As one of the most standard knot (link) invariants, Ghys investigated the linking number between the modular knot $C_\gamma$ and the ``missing" trefoil $K$ with some orientation in the $3$-sphere $S^3$. He showed that the linking number $\mathrm{Lk}(C_\gamma, K)$ is equal to the Rademacher symbol $\Psi(\gamma)$, and provided three different proofs for this result (one of them is reviewed in~\cite[Appendix A]{DIT17} and \cite{MatsusakaUeki2023}). Moreover, since there are many explicit formulas for $\Psi(\gamma)$ as above, we can compute the linking number concretely. At the end of his article, Ghys posed a question on linking numbers $\mathrm{Lk}(C_\gamma, C_\sigma)$ for two hyperbolic elements $\gamma, \sigma \in \Gamma$. In this context, a hyperbolic analogue of the Rademacher symbol $\Psi_\gamma(\sigma)$ was introduced by Duke--Imamo\={g}lu--T\'{o}th~\cite{DIT17} in 2017.

The main purpose of this article is to study the hyperbolic analogue. For two hyperbolic elements $\gamma, \sigma$, the symbol $\Psi_\gamma(\sigma)$ is defined. The details of Duke--Imamo\={g}lu--T\'{o}th's construction will be reviewed in \cref{s4-2}, but their crucial idea is to construct the weight $0$ cocycle $R_\gamma(\sigma, z)$ as a hyperbolic analogue of the cocycle $R(\sigma, z) := \log \Delta(\sigma z) - \log \Delta(z)$. Then the hyperbolic Dedekind and Rademacher symbols are defined by
\[
	\Phi_\gamma(\sigma) := \frac{2}{\pi} \lim_{y \to \infty} \Im R_\gamma(\sigma, iy), \quad \Psi_\gamma(\sigma) := \lim_{n \to \infty} \frac{\Phi_\gamma(\sigma^n)}{n},
\]
respectively. These are naturally analogous to the classical formulas
\[
	\Phi(\sigma) = \frac{1}{2\pi} \lim_{y \to \infty} \Im R(\sigma, iy), \quad \Psi(\sigma) = \lim_{n \to \infty} \frac{\Phi(\sigma^n)}{n}
\]
for a hyperbolic $\sigma \in \Gamma$. Unfortunately, the $\Psi_\gamma(\sigma)$ does not answer Ghys' original question but satisfies a similar property on linking numbers. To make this precise, we need to define two types of knots $C_\gamma^+$ and $C_\gamma^-$, whose sum is null-homologous in $\Gamma \backslash \SL_2(\mathbb{R})$~\cite[Lemma 6.1]{DIT17}. Then the hyperbolic Rademacher symbol satisfies $\Psi_\gamma(\sigma) = \mathrm{Lk}(C_\gamma^+ + C_\gamma^-, C_\sigma^+ + C_\sigma^-)$,~\cite[Theorem 3]{DIT17}. Rickards~\cite{Ric19} also gave a further study in this context. 

Get back to the classical arithmetic results, we aim to establish some expressions of $\Psi_\gamma(\sigma)$ corresponding to \eqref{Ded-Rade-rel}, \eqref{Eint-Psi}, and \eqref{Psi-explicit}. First, we get a hyperbolic analogue of the mock modular form $E_2(z)$. We recall that the $E_2(z)$ is a part of the limit of the (parabolic) Eisenstein series,
\[
	\lim_{s \to 0^+} \frac{1}{2} \sum_{\substack{c,d \in \mathbb{Z} \\ (c,d)=1}} \frac{y^s}{(cz+d)^2 \lvert cz+d \rvert^{2s}} = 1 - 24 \sum_{m=1}^\infty \left(\sum_{d \mid m} d \right) q^m - \frac{3}{\pi y} =: E_2^*(z).
\]
In \cref{s3}, we introduce a hyperbolic analogue of the above Eisenstein series, and we give a similar limit formula as follows. 

\begin{thm}[{\cref{limit-formula-for-hyperbolic}}]\label{main3}
	For a primitive hyperbolic element $\gamma = \smat{a & b \\ c & d} \in \Gamma$ with $c > 0$ and $a+d > 2$, we define the hyperbolic Eisenstein series $E_\gamma(z,s)$. Then the limit of $E_\gamma(z,s)$ as $s \to 0^+$ exists, and we have
	\[
		\lim_{s \to 0^+} E_\gamma(z,s) = -\frac{2}{\sqrt{D_\gamma}} \bigg(F_\gamma(z) - 2\log \xi_\gamma \cdot E_2^*(z) \bigg)
	\]
	for some holomorphic function $F_\gamma(z)$. The definitions of some notations are given in \cref{s3-3}. This is the $i\infty$-component of a locally harmonic Maass form of weight $2$ on $\Gamma$.
\end{thm}

The function $F_\gamma(z)$ coincides with the Duke--Imamo\={g}lu--T\'{o}th weight $2$ modular integral with rational period function studied in~\cite{DIT11, DIT17}, which is the generating series of cycle integrals of the modular functions $j_m$. We now regard the function $(2\pi i)^{-1} F_\gamma(z)$ as a hyperbolic analogue of $E_2(z)$. Then we get the following integral expression of $\Psi_\gamma(\sigma)$.

\begin{thm}[{\cref{main-theorem4}}] \label{main4}
	Let $\gamma = \smat{a & b \\ c & d} \in \Gamma$ be a primitive hyperbolic element such that $c > 0$ and $a+d > 2$. Then for a hyperbolic element $\sigma$ and $z_0 \in \bbH$, we have
	\[
		\Psi_\gamma(\sigma) = 4 \lim_{n \to \infty} \Re \int_{\sigma^n z_0}^{\sigma^{n+1} z_0} \frac{1}{2\pi i} F_\gamma(z) dz.
	\]
\end{thm}

Furthermore, the following two explicit formulas enable us to compute the hyperbolic Rademacher symbols, that is, the linking numbers $\mathrm{Lk}(C_\gamma^+ + C_\gamma^-, C_\sigma^+ + C_\sigma^-)$.

\begin{thm}[{\cref{main-theorem1}}] \label{Main-theorem1}
	Let $\gamma = \smat{a & b \\ c & d}, \sigma \in \Gamma$ be primitive hyperbolic elements such that $c > 0, a+d > 2$. We put $w_\sigma^\infty := \lim_{n \to \infty} \sigma^n z \in \{w_\sigma, w'_\sigma\}$ for any $z \in \bbH$. Then we have
	\begin{align*}
		\Psi_\gamma(\sigma) &= \Phi_\gamma(\sigma) + 2 \sum_{\substack{g \in \Gamma_{w_\gamma} \backslash \Gamma \\ w'_{g^{-1} \gamma g} < \sigma^{-1} i\infty, w_\sigma^\infty < w_{g^{-1} \gamma g}}} 1,
	\end{align*}
	where $\Gamma_{w_\gamma}$ is the stabilizer subgroup. In particular, the sum is a finite sum.
\end{thm}

This theorem says that the difference $\Psi_\gamma(\sigma) - \Phi_\gamma(\sigma)$ counts twice the number of geodesics $S_{g^{-1} \gamma g}$ intersecting with two vertical lines $x = \sigma^{-1} i\infty$ and $x = w_\sigma^\infty$. Since the hyperbolic Dedekind symbol $\Phi_\gamma(\sigma)$ also counts the number of geodesics (\cref{DIT-Ded}), we get an explicit formula for $\Psi_\gamma(\sigma)$ in terms of the counting functions. As we remark after \cref{main-theorem1}, the counting function is computable using Mathematica, for example.

\begin{thm}[{\cref{explicit-psi}}]\label{Main-theorem2}
	Let $\gamma, \sigma$ be of the form
	\[
		\gamma = \pmat{a_0 & 1 \\ 1 & 0} \cdots \pmat{a_{2n-1} & 1 \\ 1 & 0}, \qquad \sigma = \pmat{b_0 & 1 \\ 1 & 0} \cdots \pmat{b_{2m-1} & 1 \\ 1 & 0}
	\]
	with positive integers $a_i, b_j$. Then there exists an explicit function $\psi_\gamma(\sigma) \in \mathbb{Z}$ with $0 \leq \psi_\gamma(\sigma) \leq 2mn$ such that
	\[
		\Psi_\gamma(\sigma) = -2 \left(\sum_{\substack{0 \leq i < 2n \\ 0 \leq j < 2m}} \min(a_i, b_j) - \psi_\gamma(\sigma) \right).
	\]
	In particular, the value of $\Psi_\gamma(\sigma)$ is always a negative even integer.
\end{thm}

The extra term $\psi_\gamma(\sigma)$ indicates that the relationship between two modular knots is much more complicated than between a modular knot and the missing trefoil.

Finally, we mention some recent related progressions. L\"{a}geler--Schwagenscheidt~\cite{LagelerSchwagenscheidt2022} introduced a higher-weight analogue of the hyperbolic Rademacher symbol by using ``homogenized" cycle integrals considered in \cref{main4}. For studies of modular knots, Simon~\cite{Simon2022} gave further progress on Ghys' question on an arithmetic interpretation of the linking numbers of two modular knots. Ueki and the author~\cite{MatsusakaUeki2023} generalized Ghys' results for the trefoil knot and the classical Rademacher symbol to the knot complements of any torus knot using harmonic Maass forms for the triangle groups.

The rest of the article is organized as follows. In the next section, we study the Eisenstein series of weight $2$ associated with a non-scalar element $\gamma \in \Gamma$. The Eisenstein series are classified into three cases according to types of $\gamma \in \Gamma$, that is, parabolic, elliptic, and hyperbolic Eisenstein series. In particular, the limit formulas for the three types of Eisenstein series $E_\gamma(z,s)$ give contrastively a harmonic, polar harmonic, and locally harmonic Maass form. In \cref{s4}, we first review Duke--Imamo\={g}lu--T\'{o}th's construction of the hyperbolic Rademacher symbol $\Psi_\gamma(\sigma)$. After that, we establish another characterization of $\Psi_\gamma(\sigma)$ in \cref{Psi-R-red}. As corollaries, we obtain \cref{main4} and \cref{Main-theorem1}. In \cref{s5}, we provide the complete version of \cref{Main-theorem2} and its proof by using the recurrence formula for the hyperbolic Dedekind symbols. As appendices, we review the theory of continued fractions in \cref{s2} and Hooley's results on a certain exponential sum in \cref{Appendix-B}.



\section{Eisenstein series of weight $2$} \label{s3}


In this section, we study the real analytic Eisenstein series $E_\gamma(z,s)$ of weight $2$ associated with a non-scalar element $\gamma \in \Gamma := \SL_2(\mathbb{Z})$. The prototype of such a series was introduced by Petersson~\cite{Pet44}. To define it explicitly, we now prepare some notations.

Let $\mathcal{Q}$ be the set of binary quadratic forms
\[
	Q(X, Y) = [A, B, C] = A X^2 + B XY + C Y^2
\]
with integral coefficients $A, B, C \in \mathbb{Z}$. The group $\Gamma$ acts on the set $\mathcal{Q}$ by
\[
	\left(Q \circ \pmat{a & b \\ c & d} \right) (X, Y) := Q(aX + bY, cX + dY).
\]
Two quadratic forms $Q$ and $Q'$ are said to be \emph{$\Gamma$-equivalent}, and we write $Q \sim Q'$ if there exists an element $g \in \Gamma$ such that $Q' = Q \circ g$. For each element $\gamma = \smat{a & b \\ c & d} \in \Gamma$, we define the corresponding quadratic form by
\begin{align}\label{quad}
	Q_\gamma(X, Y) := cX^2 + (d-a)XY - bY^2.
\end{align}
We easily check that
\[
	Q_{g^{-1} \gamma g}(z,1) = \left( Q_\gamma \circ g \right) (z,1) = j(g, z)^2 Q_\gamma (gz, 1)
\]
for any $g \in \Gamma$. Here we put $j(g, z) := cz+d$ for any $g = \smat{a & b \\ c & d} \in \Gamma$ and $z \in \mathbb{C}$. For each quadratic form $[A, B, C] \in \mathcal{Q}$, we put $\sgn([A,B,C]) := \sgn(A)$, which takes values in $\{-1, 0, 1\}$. Under these notations, the Eisenstein series is defined as follows.

\begin{defi}\label{def-of-Eisenstein}
	Let $\pm I \neq \gamma \in \Gamma$ be a non-scalar element. For $z \in \bbH$ and $\Re(s) > 0$, the real analytic Eisenstein series is defined by
	\begin{align}\label{Def-Eisen}
	\begin{split}
		E_\gamma (z,s) &:= \sum_{Q \sim Q_\gamma} \frac{\sgn(Q) y^s}{Q(z,1) \lvert Q(z,1)\rvert^s}\\
			&= \sum_{g \in \Gamma_{w_\gamma} \backslash \Gamma} \frac{\sgn(Q_{g^{-1} \gamma g}) \Im(gz)^s}{j(g,z)^2 Q_\gamma (gz,1) \lvert Q_\gamma (gz,1) \rvert^s},
	\end{split}
	\end{align}
	where $w_\gamma$ is a fixed point of $\gamma$, and $\Gamma_{w_\gamma}$ is the stabilizer subgroup.
\end{defi}
This series converges absolutely and locally uniformly for $\Re(s) > 0$ and $z \in \bbH$. Since $Q_{-\gamma}(X,Y) = Q_{\gamma^{-1}}(X,Y) = -Q_\gamma(X,Y)$, the equation $E_{\pm \gamma^{\pm 1}} (z,s) = E_\gamma(z,s)$ holds. In addition, the series can be determined for a conjugacy class of $\gamma$, that is, $E_{g^{-1} \gamma g}(z, s) = E_\gamma(z,s)$ holds for any $g \in \Gamma$. 

For a parabolic element $\gamma$, the Eisenstein series is (essentially) the classical real analytic Eisenstein series. The elliptic case was originally investigated by Petersson. More recently, Bringmann--Kane~\cite{BK16} revisited the elliptic case in the context of polar harmonic Maass forms. As for a hyperbolic $\gamma$, there are two types of definitions. The first one is very similar to \eqref{Def-Eisen} but without the term $\sgn(Q)$. This type was studied by Zagier~\cite{Zag75-2}, Kohnen~\cite{Koh85}, and so on. Kudla--Millson~\cite{KM79} also introduced a form-valued weight $2$ analogue for each closed geodesic on the Riemann surface $\Gamma \backslash \bbH$. On the other hand, our definition \eqref{Def-Eisen} derives from Parson's work~\cite{Par93}. Parson introduced this type of hyperbolic Eisenstein series to construct a modular integral of weight $2k \geq 4$. After that, Duke--Imamo\={g}lu--T\'{o}th~\cite{DIT10} expressed its Fourier coefficients in terms of the cycle integrals of a weakly holomorphic modular form. However, they also described that it appears to be difficult to make Hecke's convergence trick work for weight $2k = 2$.

Our goal in this section is to give the limit formula of the Eisenstein series (of weight $2$) as $s \to 0^+$. Now we divide the study into three cases, (i) parabolic, (ii) elliptic, and (iii) hyperbolic, and the results are stated contrastively in \cref{PARABOLIC}, \cref{ELLIPTIC}, and \cref{limit-formula-for-hyperbolic}.


\subsection{Parabolic case} \label{s3-1}


The parabolic case is classical and well-known. Good references for the classical theory of modular forms are~\cite{Iwa97,Miy06,Zag08}, and so on. We recall that an element $\gamma \in \Gamma$ is said to be parabolic if $\gamma \neq \pm I$ and $\lvert \tr(\gamma) \rvert = 2$. The parabolic element $\gamma$ has the unique fixed point $w_\gamma$ in $\mathbb{Q} \cup \{i \infty\}$, and we can take a scaling matrix $M_\gamma \in \Gamma$ such that $M_\gamma i\infty = w_\gamma$. Since the stabilizer subgroup $\Gamma_\infty := \{g \in \Gamma \mid g i\infty = i\infty\}$ is given by $\Gamma_\infty = \{ \pm T^n \mid n \in \mathbb{Z}\}$, we have the expression $\gamma = M_\gamma T^n M_\gamma^{-1}$ or $-M_\gamma T^n M_\gamma^{-1}$ for some $n \in \mathbb{Z}_{\neq 0}$, where $T := \smat{1 & 1 \\ 0 & 1}$. By these facts, $E_\gamma(z,s) = E_{T^{\lvert n \rvert}}(z,s)$ holds, so our problem is reduced to the case when $\gamma = T^n$ with a positive integer $n > 0$. Then, the Eisenstein series is expressed as
\[
	E_{T^n}(z,s) = \sum_{g \in \Gamma_\infty \backslash \Gamma} \frac{\sgn(Q_{g^{-1} T^n g}) y^s}{Q_{g^{-1} T^n g} (z,1) \lvert Q_{g^{-1} T^n g} (z,1) \rvert^s}.
\]
For each $g = \smat{* & * \\ c & d} \in \Gamma_\infty \backslash \Gamma$, we see that $g^{-1} T^n g = \smat{1+cdn & d^2 n \\ -c^2n & 1-cdn}$, that is, $Q_{g^{-1} T^n g}(z,1) = -n(cz+d)^2 = -n j(g,z)^2$. It implies the expression
\[
	E_{T^n}(z,s) = \frac{1}{n^{s+1}} \sum_{\substack{g \in \Gamma_\infty \backslash \Gamma \\ g \neq I}} \frac{\Im(gz)^s}{j(g,z)^2} = \frac{1}{n^{s+1}} E_T(z,s).
\]
Here we remark that $\sgn(Q_{g^{-1} T^n g}) = 0$ holds for any $g \in \Gamma_\infty$. By using Hecke's trick or the Fourier expansion, we get the analytic continuation of the function to $s = 0$. 
\begin{prop}~\cite[Proposition\ 6]{Zag08} \label{PARABOLIC}
	For a parabolic element $T = \smat{1 & 1 \\ 0 & 1}$, we have
	\[
		\lim_{s \to 0^+} E_T(z,s) = - 24 \sum_{m=1}^\infty \sigma_1(m) q^m - \frac{3}{\pi y},
	\]
	where $\sigma_1(m) := \sum_{d \mid m} d$ is the divisor sum. 
\end{prop}

The series $E_2^*(z) = 1 + \lim_{s \to 0^+} E_T(z,s)$ is a harmonic Maass form of weight $2$ on $\Gamma$. Here a real analytic function $f: \bbH \to \mathbb{C}$ is called a \emph{harmonic Maass form} of weight $k \in \mathbb{Z}$ on $\Gamma$ if
\begin{enumerate}
	\item $j(g,z)^{-k} f(gz) = f(z)$ for any $g \in \Gamma$,
	\item $\displaystyle{\Delta_k f(z) := \left(-y^2 \left(\frac{\partial^2}{\partial x^2} + \frac{\partial^2}{\partial y^2} \right) + iky \left(\frac{\partial}{\partial x} + i \frac{\partial}{\partial y} \right) \right) f(z) = 0}$,
	\item $f(x+iy) = O(y^\alpha)$ as $y \to \infty$ for some $\alpha > 0$.
\end{enumerate}


\subsection{Elliptic case} \label{s3-2}


An element $\gamma \in \Gamma$ is said to be elliptic if $\lvert \tr(\gamma)\rvert < 2$. The elliptic element $\gamma$ has two fixed points $w_\gamma \in \bbH$ and its complex conjugate $\overline{w}_\gamma$. The remarkable fact is that any elliptic fixed point $w_\gamma$ is $\Gamma$-equivalent to either $i$ or $\rho := e^{\pi i/3}$, that is, $g w_\gamma = i$ or $\rho$ holds for some $g \in \Gamma$. Moreover the stabilizer subgroups $\Gamma_{w_\gamma}$ for $w_\gamma = i, \rho$ are finite cyclic groups given by $\Gamma_i = \{S^n := \smat{0 & -1 \\ 1 & 0}^n \mid 0 \leq n \leq 3\}$ and $\Gamma_\rho = \{U^n := \smat{1 & -1 \\ 1 & 0}^n \mid 0 \leq n \leq 5\}$. Therefore, in this case, our problem is reduced to the cases of $\gamma = S$ and $U$.

For $\gamma = S$, the Eisenstein series is expressed as
\begin{align*}
	E_S(z,s) &= \frac{1}{\lvert \Gamma_i \rvert} \sum_{g \in \Gamma} \frac{\sgn(Q_{g^{-1} S g}) \Im(gz)^s}{j(g,z)^2 Q_S (gz,1) \lvert Q_S (gz,1) \rvert^s}\\
		&= \frac{1}{4} \sum_{g \in \Gamma} \frac{\Im(gz)^s}{j(g,z)^2 (gz-i)(gz+i) \lvert (gz-i)(gz+i) \rvert^s}.
\end{align*}
Here $\sgn(Q_{g^{-1} S g})$ is always positive. In general, for any $z, \tau \in \bbH$, which are $\Gamma$-inequivalent each other, we introduce another Poincar\'{e} series
\begin{align}\label{E-elliptic}
	E(z, \tau; s) := \sum_{g \in \Gamma} \frac{\Im(\tau)^{s+1} \Im(gz)^s}{j(g,z)^2 (gz-\tau)(gz-\overline{\tau}) \lvert (gz-\tau)(gz-\overline{\tau})\rvert^s}
\end{align}
for $\Re(s) > 0$. Then our Eisenstein series is expressed as $E_\gamma(z,s) = \Im(w_\gamma)^{-s-1} E(z, w_\gamma; s) / \lvert \Gamma_{w_\gamma}\rvert$ for $\gamma = S, U$. The function $E(z, \tau; s)$ satisfies the modular transformation laws $j(g,z)^{-2} E(gz, \tau; s) = E(z, \tau; s)$ and $E(z, g\tau; s) = E(z, \tau; s)$ for any $g \in \Gamma$. These immediately follow from the facts that $j(gh, z) = j(g, hz) j(h,z)$ and $j(g, \tau) (z - g\tau) = j(g^{-1}, z) (g^{-1}z - \tau)$.

Now we give the analytic continuation of $E(z, \tau; s)$ to $s = 0$ generally, which was essentially established by Bringmann--Kane~\cite{BK16} and others~\cite{BKLOR18} as follows. First, Bringmann--Kane considered a slightly different Poincar\'{e} series
\[
	\mathcal{P}(z, \tau; s) := \sum_{g \in \Gamma} \frac{\Im(\tau)^{s+1}}{j(g, z)^2 \lvert j(g,z)\rvert^{2s} (gz - \tau) (gz - \bar{\tau}) \lvert gz - \bar{\tau} \rvert^{2s}}
\]
for $\Re(s) > 0$. After a lengthy calculation, they derived the analytic continuation of $\mathcal{P}(z, \tau; s)$ to $s = 0$, and the following explicit formula.
\begin{align}\label{BK-limit}
	\mathcal{P}(z, \tau; 0) = -2\pi \left(\frac{j'(z)}{j(\tau) - j(z)} - E_2^*(z) \right),
\end{align}
where $j(z)$ is the elliptic modular $j$-function defined by
\[
	j(z) := \frac{\displaystyle{\left(1 + 240 \sum_{d=1}^\infty \frac{d^3 q^d}{1-q^d} \right)^3}}{\displaystyle{q \prod_{n=1}^\infty (1- q^n)^{24}}} = q^{-1} + 744 + 196884q + \cdots,
\]
and $' = (2\pi i)^{-1} d/dz$. The function \eqref{BK-limit} in $z$ gives a polar harmonic Maass form of weight $2$ on $\Gamma$, which is a harmonic Maass form allowed to have poles in the upper half-plane. For more details, see~\cite[Section 13.3]{BFOR17}. 

Next, we compare these two Poincar\'{e} series in $\Re(s) > 0$,
\begin{align}\label{E-P}
	\begin{split}
		&E(z, \tau; s) - \Im(z)^s \mathcal{P}(z, \tau; s)\\
			&= \sum_{g \in \Gamma} \frac{\Im(\tau)^{s+1} \Im(z)^s}{j(g, z)^2 \lvert j(g,z) \rvert^{2s} (gz - \tau) (gz - \bar{\tau}) \lvert gz - \bar{\tau}\rvert^{2s}} \left( \left\lvert \frac{gz - \tau}{gz - \bar{\tau}} \right\rvert^{-s} - 1 \right).
	\end{split}
\end{align}
To estimate the last term, we need the following lemma.

\begin{lem}\label{evaluate}
	For a fixed $\tau \in \bbH$, let $K \subset \bbH \setminus \{g\tau \in \bbH \mid g \in \Gamma\}$ be any compact subset. Then, there exists a constant $C = C(\tau, K) > 0$ such that, for any $s \in [-1,1]$ and $z \in K$,
	\[
		\left\lvert \left\lvert \frac{gz - \tau}{gz - \bar{\tau}} \right\rvert^{-s} - 1 \right\rvert \leq C \cdot \frac{\Im(gz)}{\lvert gz - \bar{\tau} \rvert^2}.
	\]
\end{lem}

\begin{proof}
	We put
	\[
		\left\lvert \frac{z - \tau}{z - \bar{\tau}} \right\rvert^2 = 1 - \frac{4 \Im(z) \Im(\tau)}{\lvert z - \bar{\tau}\rvert^2} =: 1 - a_\tau(z).
	\]
	For any $\varepsilon > 0$ and a fixed $\tau$, the region $\{z \in \bbH \mid a_\tau(z) \geq \varepsilon\}$ is a compact subset in $\bbH$. Then, there exists $M = M(\tau, K) \in (0,1)$ such that $a_\tau(gz) < M$ for any $g \in \Gamma$ and $z \in K$ since $\Gamma$ acts properly discontinuously on $\bbH$. By applying the mean value theorem to $f(r) = (1-r)^{-s/2}$, we have
	\[
		\left\lvert \frac{gz - \tau}{gz - \bar{\tau}} \right\rvert^{-s} = (1 - a_\tau(gz))^{-s/2} = 1 + \frac{s}{2} (1 - c)^{-s/2-1} a_\tau(gz)
	\]
	for some $0 < c < a_\tau(gz) < M$. Therefore, for any $s \in [-1, 1]$, we have
	\[
		\left\lvert \left\lvert \frac{gz - \tau}{gz - \bar{\tau}} \right\rvert^{-s} - 1 \right\rvert \leq (1-M)^{-3/2} a(gz).
	\]
	This concludes the proof.
\end{proof}

By this lemma, the right-hand side of \eqref{E-P} converges absolutely and locally uniformly for $s \in (-1,1]$, and has a zero at $s = 0$. Thus we get the following results.

\begin{thm}\label{limit-E}
	Let $z, \tau \in \bbH$ be $\Gamma$-inequivalent elements. Then we have
	\[
		\lim_{s \to 0^+} E(z, \tau; s) = -2\pi \left( \frac{j'(z)}{j(\tau) - j(z)} - E_2^*(z) \right),
	\]
	which gives a polar harmonic Maass form of weight $2$ on $\Gamma$ in $z$, whose poles are located on the $\Gamma$-orbit of $\tau$.
\end{thm}

\begin{cor}\label{ELLIPTIC}
	For elliptic elements $S = \smat{0 & -1 \\ 1 & 0}$ and $U = \smat{1 & -1 \\ 1 & 0}$, we have
\begin{align*}
	\lim_{s \to 0^+} E_S(z,s) &= \frac{\pi}{2} \left( \frac{j'(z)}{j(z)-1728} + E_2^*(z) \right),\\
	\lim_{s \to 0^+} E_U(z,s) &= \frac{2\pi}{3\sqrt{3}} \left( \frac{j'(z)}{j(z)} + E_2^*(z) \right),
\end{align*}
	which give polar harmonic Maass forms of weight $2$ on $\Gamma$.
\end{cor}
Finally, we note that Asai--Kaneko--Ninomiya~\cite{AKN97} showed the Fourier expansion
\begin{align}\label{AKN}
	\frac{j'(z)}{j(\tau) - j(z)} = \sum_{m=0}^\infty j_m(\tau) q^m, \qquad \Im(z) > \Im(\tau),
\end{align}
where $j_m(z)$ is the unique polynomial in $j(z)$ having a Fourier expansion of the form $q^{-m} + O(q)$. The first few examples are 
\begin{align*}
	j_0(z) &= 1,\\
	j_1(z) &= j(z) - 744 = q^{-1} + 196884q + \cdots,\\
	j_2(z) &= j(z)^2 - 1488j(z) + 159768 = q^{-2} + 42987520q + \cdots.
\end{align*}
The set $\{j_m(z)\}_{m=0}^\infty$ consists a natural basis for the space $M_0^!(\Gamma)$ of weakly holomorphic modular forms of weight $0$ on $\Gamma$, and form a Hecke system in the sense of~\cite{BKLOR18}. Combining with the Fourier expansion of $E_2^*(z)$, we get the expansion
\[
	\frac{j'(z)}{j(\tau) - j(z)} - E_2^*(z) = \sum_{m=1}^\infty (j_m(\tau) + 24 \sigma_1(m)) q^m + \frac{3}{\pi y}
\]
for $\Im(z) > \Im(\tau)$.


\subsection{Hyperbolic case} \label{s3-3}


This case is slightly different from the above two cases. An element $\gamma \in \Gamma$ is said to be hyperbolic if $\lvert \tr(\gamma) \rvert > 2$. The hyperbolic element $\gamma$ has two fixed points $w_\gamma, w'_\gamma$ which are real quadratic irrationals and Galois conjugate each other. In this article, we always order them by $w_\gamma > w'_\gamma$. We call $\gamma$ \emph{primitive} if $\Gamma_{w_\gamma} :=\{g \in \Gamma \mid g w_\gamma = w_\gamma\} = \{ \pm \gamma^n \mid n \in \mathbb{Z}\}$. Furthermore, we denote by $S_\gamma$ the oriented geodesic in $\bbH$ connecting two fixed points $w_\gamma$ and $w'_\gamma$. Here the orientation is clockwise if $\sgn(\gamma) := \sgn(c(a+d)) = \sgn(Q_\gamma) \sgn(\mathrm{tr}(\gamma)) > 0$, and counter-clockwise if $\sgn(\gamma) < 0$. Unfortunately, the Eisenstein series $E_\gamma(z,s)$ does not satisfy the modular transformation law. Actually, the series $E_\gamma(z,s)$ satisfies the following.

\begin{prop}\label{Egamma-gap}
	Let $\gamma$ be a hyperbolic element. Then for $\Re(s) > 0$ and any $\sigma = \smat{a & b \\ c & d} \in \Gamma \setminus \Gamma_\infty$, we have
	\begin{align}\label{Ehyp-gap}
		j(\sigma, z)^{-2} E_\gamma(\sigma z, s) - E_\gamma(z,s) = -2 \sum_{\substack{g \in \Gamma_{w_\gamma} \backslash \Gamma \\ w'_{g^{-1} \gamma g} < \sigma^{-1} i\infty < w_{g^{-1} \gamma g}}} \frac{\sgn(Q_{g^{-1} \gamma g}) \Im(z)^s}{Q_{g^{-1} \gamma g}(z,1) \lvert Q_{g^{-1} \gamma g}(z,1) \rvert^s},
	\end{align}
	where $\sigma^{-1} i\infty = -d/c \in \mathbb{Q}$. The sum in the right-hand side is a finite sum. For $\sigma \in \Gamma_\infty$, the right-hand side equals $0$.
\end{prop}

\begin{proof}
	The finiteness is a known fact. We will describe the details after \cref{main-theorem1}. By the definition of the Eisenstein series,
	\begin{align*}
		j(\sigma, z)^{-2} E_\gamma(\sigma z, s) &= \sum_{g \in \Gamma_{w_\gamma} \backslash \Gamma} \frac{\sgn(Q_{g^{-1} \gamma g}) \Im(z)^s}{Q_{\sigma^{-1} g^{-1} \gamma g \sigma}(z,1) \lvert Q_{\sigma^{-1} g^{-1} \gamma g \sigma}(z,1) \rvert^s}\\
			&= \sum_{g \in \Gamma_{w_\gamma} \backslash \Gamma} \frac{\sgn(Q_{\sigma g^{-1} \gamma g \sigma^{-1}}) \Im(z)^s}{Q_{g^{-1} \gamma g}(z,1) \lvert Q_{g^{-1} \gamma g}(z,1)\rvert^s}.
	\end{align*}
	Thus 
	\[
		j(\sigma, z)^{-2} E_\gamma (\sigma z, s) - E_\gamma(z,s) = \sum_{g \in \Gamma_{w_\gamma} \backslash \Gamma} \frac{(\sgn(Q_{\sigma g^{-1} \gamma g \sigma^{-1}}) - \sgn(Q_{g^{-1} \gamma g})) \Im(z)^s}{Q_{g^{-1} \gamma g}(z,1) \lvert Q_{g^{-1} \gamma g}(z,1)\rvert^s}.
	\]
	For $\sigma \in \Gamma_\infty$, we easily see that $\sgn(Q_{\sigma g^{-1} \gamma g \sigma^{-1}}) = \sgn(Q_{g^{-1} \gamma g})$, that is, the above equals $0$. We now show that
	\begin{align*}
		\sgn(Q_{\sigma h \sigma^{-1}}) \sgn(Q_h) = \begin{cases}
			-1 \qquad &\text{if } w'_h < \sigma^{-1} i\infty < w_h,\\
			+1 &\text{if otherwise}
		\end{cases}
	\end{align*}
	for any $\sigma \in \Gamma \setminus \Gamma_\infty$ and any hyperbolic $h$. This immediately implies \eqref{Ehyp-gap}. Let $h = \smat{A & B \\ C & D} \in \Gamma$ be hyperbolic. By a direct calculation, we have
	\begin{align*}
		\sgn(Q_{\sigma h \sigma^{-1}}) \sgn(Q_h) &= \sgn((Acd - Bc^2 + Cd^2 - Dcd)C)\\
			&= \sgn((Ac-Dc+2Cd)^2 - c^2 ((A+D)^2-4)).
	\end{align*}
	Thus, the product is negative if and only if
	\[
		\lvert Ac - Dc + 2Cd \rvert < \lvert c\rvert \sqrt{(A+D)^2 - 4},
	\]
	that is, $w'_h < \sigma^{-1} i\infty < w_h$.
\end{proof}

We next give the limit formula for the Eisenstein series $E_\gamma(z,s)$ as $s \to 0^+$. To describe the result, we introduce more notations. First, we show the following lemmas.

\begin{lem}\label{gamma-limit}
	Let $\gamma = \smat{a & b \\ c & d} \in \Gamma$ be a hyperbolic element with fixed points $w_\gamma > w'_\gamma$. For any $z \in \bbH \cup \mathbb{R} \cup \{i\infty\} \setminus \{w_\gamma, w'_\gamma\}$, we have
	\[
		\lim_{n \to \infty} \gamma^n z = \begin{cases}
			w_\gamma &\text{if } \sgn(\gamma) := \sgn(c(a+d)) > 0,\\
			w'_\gamma & \text{if } \sgn(\gamma) < 0.
		\end{cases}
	\]
	It is compatible with the orientation of the geodesic $S_\gamma$.
\end{lem} 

\begin{proof}
	Since $\gamma$ is diagonalized by the scaling matrix $M_\gamma$ as in \eqref{fund-unit}, we have
	\[
		\gamma^n = M_\gamma \pmat{\xi_\gamma^n & 0 \\ 0 & \xi_\gamma^{-n}} M_\gamma^{-1} = \frac{1}{w_\gamma - w'_\gamma} \pmat{w_\gamma \xi_\gamma^n - w'_\gamma \xi_\gamma^{-n} & w_\gamma w'_\gamma (\xi_\gamma^{-n} - \xi_\gamma^n) \\ \xi_\gamma^n - \xi_\gamma^{-n} & -w'_\gamma \xi_\gamma^n + w_\gamma \xi_\gamma^{-n}}.
	\]
	Here $\xi_\gamma := cw_\gamma + d = \dfrac{a+d + \sgn(c) \sqrt{(a+d)^2 - 4}}{2}$ satisfies $\lvert \xi_\gamma\rvert > 1$ if $c(a+d) > 0$, and $0 < \lvert \xi_\gamma \rvert < 1$ if $c(a+d) < 0$.
\end{proof}

\begin{lem}
	For a hyperbolic element $\gamma \in \Gamma$, we suppose that $\sgn(Q_\gamma) > 0$ and $\mathrm{tr}(\gamma) > 2$. For a positive integer $n > 0$, we have
	\[
		E_{\gamma^n} (z,s) = \frac{1}{a_n^{s+1}} E_\gamma(z,s),
	\]
	where the sequence $\{a_n\}_{n=1}^\infty$ is defined by $a_n = \mathrm{tr}(\gamma) a_{n-1} - a_{n-2}$ with the initial values $a_0 = 0, a_1 = 1$.
\end{lem}

\begin{proof}
	By the definition \eqref{Def-Eisen} of the Eisenstein series, 
	\[
		E_{\gamma^n}(z,s) = \sum_{g \in \Gamma_{w_{\gamma^n}} \backslash \Gamma} \frac{\sgn(Q_{g^{-1} \gamma^n g}) \Im(gz)^s}{j(g,z)^2 Q_{\gamma^n}(gz,1) \lvert Q_{\gamma^n}(gz,1)\rvert^s}.
	\]
	We immediately see that $\Gamma_{w_{\gamma^n}} = \Gamma_{w_\gamma}$. Let $C_n(x) \in \mathbb{Z}[x]$ denote the Chebyshev polynomial of the second kind characterized by $C_0(x) = 0, C_1(x) = 1$, and $C_{n+1}(x) = 2x C_n(x) - C_{n-1}(x)$. Then, by induction, it immediately follows that
	\[
		\pmat{a & b \\ c & d}^n = \pmat{* & * \\ c C_n \left(\frac{a+d}{2}\right) & d C_n \left(\frac{a+d}{2}\right) - C_{n-1} \left(\frac{a+d}{2}\right)}.
	\]
	Since $C_n(\cos t) = \sin nt/\sin t$ holds, all zeros of $C_n(x)$ are real and located in $(-1,1)$. It implies that $C_n((a+d)/2) > 0$ for $a+d > 2$. Thus, we have $\sgn(Q_{g^{-1} \gamma^n g}) = \sgn(Q_{g^{-1} \gamma g})$. Furthermore, since $Q_{\gamma^n}$ has the same roots as $Q_\gamma$, we have
	\[
		Q_{\gamma^n} = C_n \left(\frac{a+d}{2} \right) Q_\gamma = a_n Q_\gamma,
	\]
	which completes the proof.
\end{proof}

From now on, we assume that $\gamma \in \Gamma$ is a primitive hyperbolic element with $\sgn(Q_\gamma) > 0$ and $\mathrm{tr}(\gamma) > 2$. Let $j_m(z)$ be the modular function of the form $j_m(z) = q^{-m} + O(q)$ defined in \eqref{AKN}. We now consider the cycle integral
\[
	\widetilde{\mathrm{val}}_m(\gamma) := \int_{\tau_0}^{\gamma \tau_0} j_m(\tau) \frac{-\sqrt{D_\gamma} d\tau}{Q_\gamma(\tau,1)} = \int_{\Gamma_{w_\gamma} \backslash S_\gamma}  j_m(\tau) \frac{\lvert d\tau \rvert}{\Im(\tau)},
\]
where $D_\gamma := \mathrm{tr}(\gamma)^2 - 4$ is the discriminant of $Q_\gamma(X,Y)$, and the path of integration is on the geodesic $S_\gamma$. The notation ``$\mathrm{val}$" was introduced by Kaneko~\cite{Kan09} to represent the normalized cycle integral $\mathrm{val}(w_\gamma) = \widetilde{\mathrm{val}}_1(\gamma)/\widetilde{\mathrm{val}}_0(\gamma) + 744$. In particular, $\widetilde{\mathrm{val}}_0(\gamma) = 2\log \xi_\gamma$ is the length of the closed geodesic $\Gamma_{w_\gamma} \backslash S_\gamma$. The cycle integrals of the modular $j$-function have been studied by Kaneko and Duke--Imamo\={g}lu--T\'{o}th~\cite{DIT11} as a real quadratic analogue of singular moduli. Kaneko made several observations on the distribution of $\mathrm{val}(w)$ through numerical experiments, some of which were proved by Bengoechea--Imamo\={g}lu~\cite{BI19,BI20} and Murakami~\cite{Mur19}. 

As shown in~\cite{DIT11,DIT17}, by using Hooley's estimate on a certain exponential sum, one can see that $\widetilde{\val}_m(\gamma)$ has at most polynomial growth as $m \to \infty$.
Then the $q$-series
\begin{align}\label{F-cocycle}
	F_\gamma(z) := \sum_{m=0}^\infty \widetilde{\mathrm{val}}_m(\gamma) q^m
\end{align}
converges on the whole $\bbH$-plane, and defines a holomorphic function. Under these notations, our theorem is described as follows.

\begin{thm}[\cref{main3}] \label{limit-formula-for-hyperbolic}
	For a primitive hyperbolic element $\gamma = \smat{a & b \\ c & d} \in \Gamma$ with $c> 0$ and $a+d  >2$, we have
	\begin{align}\label{limit-formula-for-hyperbolic-Eis}
		\lim_{s \to 0^+} E_\gamma(z,s) = -\frac{2}{\sqrt{D_\gamma}} \bigg(F_\gamma(z) - 2\log \xi_\gamma \cdot E_2^*(z) \bigg).
	\end{align}
	This is the $i\infty$-component of a locally harmonic Maass form of weight $2$ on $\Gamma$.
\end{thm}

A locally harmonic Maass form is a harmonic Maass form with jumping singularity on the net of geodesics in $\bbH$, which Bringmann--Kane--Kohnen~\cite{BKK15} originally introduced. The function defined by
\[
	-\frac{2}{\sqrt{D_\gamma}} \int_{\tau_0}^{\gamma \tau_0} \left(\frac{j'(z)}{j(\tau) - j(z)} - E_2^*(z) \right)\frac{-\sqrt{D_\gamma}d\tau}{Q_\gamma(\tau,1)}
\]
is a modular form of weight $2$ on $\Gamma$, annihilated by the action of $\Delta_2$, and has jumping singularities on geodesics. The function is a locally harmonic Maass form of weight $2$, and coincides with \eqref{limit-formula-for-hyperbolic-Eis} for a large enough $\Im(z)$ by \eqref{AKN}. We call the part of a locally harmonic Maass form on the connected component, including $i\infty$ without singularities, its \emph{$i\infty$-component}. For more details, see~\cite[Section 13.4]{BFOR17}. Recently, L\"{o}brich--Schwagenscheidt~\cite[Theorem 4.2]{LS20} constructed similar locally harmonic Maass forms of weight $2 - 2k \leq 0$.

\begin{rem}
	Mono~\cite{Mono2022, Mono2021pre} showed the traced-version of \cref{limit-formula-for-hyperbolic} inspired by the first version of this article.
\end{rem}

To prove \cref{limit-formula-for-hyperbolic}, we recall some results by Duke--Imamo\={g}lu--T\'{o}th~\cite{DIT11}. For fixed $\nu \in \mathbb{C}$ and $y > 0$, the Bessel functions are defined by
\begin{align}\label{IJ-Bessel}
	I_\nu(y) = \sum_{k=0}^\infty \frac{(y/2)^{\nu + 2k}}{k! \Gamma(\nu + k + 1)}, \quad J_\nu(y) = \sum_{k=0}^\infty \frac{(-1)^k (y/2)^{\nu+2k}}{k! \Gamma(\nu+k+1)}.
\end{align}
Let $m \in \mathbb{Z}$. The Niebur Poincar\'{e} series is defined by
\[
	G_m(\tau, s) := \sum_{\gamma \in \Gamma_\infty \backslash \Gamma} e^{2\pi i m \Re(\gamma \tau)} \phi_{m,s}(\Im(\gamma \tau))
\]
with
\[
	\phi_{m,s}(y) := \begin{cases}
		2\pi \lvert m\rvert^{1/2} y^{1/2} I_{s-1/2}(2\pi \lvert m \rvert y) \qquad &\text{if } m \neq 0,\\
		y^s &\text{if } m=0.
	\end{cases}
\]
This series converges absolutely and uniformly on compact subsets for $\Re(s)>1$. Moreover, the series has the meromorphic continuation to $\Re(s) > 3/4$ via its Fourier expansion. For $m > 0$, the series $G_{-m}(\tau,s)$ is holomorphic at $s =1$, and equals
\begin{align}\label{Niebur-analytic}
	G_{-m}(\tau,1) = j_m(\tau) + 24 \sigma_1(m).
\end{align}
As for $m = 0$, this series has a simple pole at $s = 1$, and satisfies
\begin{align}\label{Niebur-merom}
	G_0(\tau, s) = \frac{3}{\pi} \cdot \frac{1}{s-1} + O(1)
\end{align}
at $s = 1$. For more details, see~\cite[Section 4]{DIT11}.

For a non-zero integer $k \neq 0$ and an integer $m \in \mathbb{Z}$, we define an exponential sum
\begin{align}\label{exponential-sum}
	S_\gamma(k, m) = \sum_{\substack{\ell\ (\mathrm{mod}\ 2\lvert k \rvert) \\ [k,\ell, \frac{\ell^2 - D_\gamma}{4k}] \sim Q_\gamma}} e^{\pi i m\ell/k}.
\end{align}

\begin{lem}\label{DIT-Prop 4}
	Let $\Re(s) > 1$ and $m \in \mathbb{Z}$. Suppose that $\gamma = \smat{a & b \\ c & d} \in \Gamma$ is a primitive hyperbolic element with $c > 0$ and $a+d > 2$. Then, we have
	\begin{align*}
		&\frac{2\Gamma(s)}{\Gamma \left(\frac{s}{2}\right)^2} \int_{\tau_0}^{\gamma \tau_0} G_m(\tau, s) \frac{-\sqrt{D_\gamma} d\tau}{Q_\gamma(\tau, 1)}\\
			&= \begin{cases}
				\displaystyle{2^{s+1/2} \pi \lvert m\rvert^{1/2} D_\gamma^{1/4} \sum_{a \in \mathbb{Z}_{\neq 0}} \frac{S_\gamma(a, -m)}{\lvert a\rvert^{1/2}} J_{s-1/2} \left(\frac{\pi \lvert m \rvert \sqrt{D_\gamma}}{\lvert a\rvert} \right)} \ &\text{if } m \neq 0,\\
				\displaystyle{D_\gamma^{s/2} \sum_{a \in \mathbb{Z}_{\neq 0}} \frac{S_\gamma(a,0)}{\lvert a\rvert^s}} &\text{if } m = 0.
			\end{cases}
	\end{align*}
\end{lem}

\begin{proof}
	It follows by almost the same arguments as in the proof of~\cite[Proposition 4]{DIT11}. By a similar argument as in \cite[Lemma 7]{DIT11}, we have
	\[
		\int_{\tau_0}^{\gamma \tau_0} G_m(\tau, s) \frac{-\sqrt{D_\gamma} d\tau}{Q_\gamma(\tau,1)} = \sum_{a \in \mathbb{Z}_{\neq 0}} \Phi_m \left(s, \frac{\sqrt{D_\gamma}}{2 \lvert a\rvert} \right) S_\gamma(a, -m),
	\]
	where
	\[
		\Phi_m(s,t) = \int_0^\pi \cos(2\pi mt \cos \theta) \phi_{m,s}(t \sin \theta) \frac{d\theta}{\sin \theta}.
	\]
	The formula for the integral $\Phi_m(s,t)$ was shown in \cite[Lemma 9]{DIT11}, which gives the desired result.
\end{proof}

Furthermore, we have the evaluation of the above equation at $s=1$. The ``traced" version was shown by Duke--Imamo\={g}lu--T\'{o}th~\cite[Section 6]{DIT11}. The following ``individual" version can be shown in entirely the same way.

\begin{thm}\label{val-explicit-formula}
	Let $m > 0$ be a positive integer and $\gamma$ be as in \cref{DIT-Prop 4}. Then the series
	\[
		2^{1/2} \pi^2 m^{1/2} D_\gamma^{1/4} \sum_{a \in \mathbb{Z}_{\neq 0}} \frac{S_\gamma(a,m)}{\lvert a\rvert^{1/2}} J_{1/2} \left(\frac{\pi m \sqrt{D_\gamma}}{\lvert a\rvert} \right)
	\]
	(conditionally) converges and coincides with $\widetilde{\val}_m(\gamma) + 24 \sigma_1(m) \widetilde{\val}_0(\gamma)$. For $m = 0$, we have
	\[
		\sum_{a \in \mathbb{Z}_{\neq 0}} \frac{S_\gamma(a,0)}{\lvert a \rvert^s} = \frac{6 \widetilde{\val}_0(\gamma)}{\pi^2 \sqrt{D_\gamma}} \cdot \frac{1}{s-1} + O(1)
	\]
	as $s \to 1$.
\end{thm}

\begin{proof}
	Let $m > 0$. First, we recall Hooley's estimate~\cite{Hoo63}, for any $\epsilon > 0$,
	\begin{align}\label{Hooley}
		\sum_{0 < k \leq x} S_\gamma(k,m) \ll_\epsilon m^\epsilon x^\epsilon (m^{4/5} + x^{3/4}).
	\end{align}
	Since the proof is omitted in~\cite{Hoo63}, we review it in \cref{Appendix-B} for convenience. \footnote{In~\cite{DIT11}, there is a minor typo on the exponent of $m$ in~\eqref{Hooley}. The corrected proof gives a stronger result, $\widetilde{\val}_m(\gamma) + 24\sigma_1(m) \widetilde{\val}_0(\gamma) \ll_{\gamma, \epsilon} m^{1+\epsilon}$, than \cite[Propositoin 6]{DIT11}.}
	
	For any $s \in [1,2]$ and $t \in [3/8, 1]$, the trivial bound $S_\gamma(k,m) \ll_\epsilon \lvert k\rvert ^\epsilon$ and the bound $J_{s-1/2}(y) \leq 2y^{-1/2}$ (see Olenko~\cite{Olenko2006}) imply that
	\begin{align}\label{Est-1}
		\sum_{0 < k \leq m} S_\gamma(k,m) k^{-t} J_{s-1/2} \left(\frac{\pi m \sqrt{D_\gamma}}{k} \right) &\ll_\epsilon m^{-1/2} \sum_{0 < k \leq m} k^{\epsilon - t +1/2} \ll_\epsilon m^{1 - t + \epsilon}.
	\end{align}
	By the definition of the $J$-Bessel function, for $x > m$, 
	\begin{align}\label{Est-2}
		\sum_{m < k \leq x} S_\gamma(k,m) k^{-t} J_{s-1/2} \left(\frac{\pi m \sqrt{D_\gamma}}{k} \right) \ll_\epsilon m^{s-1/2} \left\lvert \sum_{m < k \leq x} S_\gamma(k,m) k^{-s-t+1/2} \right\rvert + m^{1-t+\epsilon}.
	\end{align}
	Abel's summation formula implies
	\begin{align*}
		\sum_{m < k \leq x} S_\gamma(k,m) k^{-s-t+1/2} &\ll x^{-s-t+1/2} \sum_{1 \leq k \leq x} S_\gamma(k,m) - m^{-s-t+1/2} \sum_{1 \leq k \leq m} S_\gamma(k,m)\\
		&\qquad + \int_m^x u^{-s-t-1/2} \sum_{1 \leq k \leq u} S_\gamma(k,m) du.
	\end{align*}
	By Hooley's estimate~\eqref{Hooley}, we have
	\[
		\sum_{m < k \leq x} S_\gamma(k,m) k^{-s-t+1/2} \ll_\epsilon m^{1/2+4/5-s-t+\epsilon}.
	\]
	Combining it with \eqref{Est-1} and \eqref{Est-2}, we obtain
	\begin{align}\label{Hooley-DIT}
	\begin{split}
		\sum_{0 < k \leq x} S_\gamma(k,m) k^{-t} J_{s-1/2} \left(\frac{\pi m \sqrt{D_\gamma}}{k} \right) &\ll_\epsilon m^{1-t+\epsilon} + m^{4/5 - t + \epsilon}\\
			&\ll_\epsilon m^{1-t+\epsilon}.
	\end{split}
	\end{align}
	In particular, it is bounded on $x > 0$. By applying Dirichlet's test to the monotonically decreasing function $k^{-1/8}$ and the above uniformly bounded sequence at $t=3/8$ on $s \in [1,2]$, the series
	\[
		\sum_{a \in \mathbb{Z}_{\neq 0}} \frac{S_\gamma(a,m)}{\lvert a \rvert^{1/2}} J_{s-1/2} \left(\frac{\pi m \sqrt{D_\gamma}}{\lvert a\rvert} \right)
	\]
	uniformly converges on $s \in [1,2]$. We conclude the results by taking the limit $s \to 1$ in \cref{DIT-Prop 4} for both cases $m>0$ and $m=0$. Here we note the relation $S_\gamma(-a,m) = S_{\gamma^{-1}}(a,m)$.
\end{proof}

\begin{proof}[Proof of \cref{limit-formula-for-hyperbolic}]
	For $\Re(s) > 0$, we split the sum in the definition of $E_\gamma(z,s)$ as follows.
	\begin{align}\label{Egamma-decomposition}
	\begin{split}
		E_\gamma(z,s) &= \sum_{[a,b,c] \sim Q_\gamma} \frac{\sgn(a) y^s}{(az^2+bz+c)\lvert az^2 + bz + c\rvert^s}\\
		&= y^s \sum_{a \in \mathbb{Z}_{\neq 0}} \frac{1}{\lvert a \rvert^{s+1}} \sum_{\substack{b\ (2\lvert a\rvert) \\ [a,b,\frac{b^2-D_\gamma}{4a}] \sim Q_\gamma}}  \sum_{n \in \mathbb{Z}} \frac{1}{\left(z+n + \frac{b}{2a}\right)^2 - \frac{D_\gamma}{4a^2}} \frac{1}{\left\lvert \left(z+n + \frac{b}{2a}\right)^2 - \frac{D_\gamma}{4a^2} \right\rvert^s}.
	\end{split}
	\end{align}
	We consider the most inner summand. Let $w \in \mathbb{H}$ with fixed imaginary part $y > 0$, and $r \in \mathbb{R}$. For $s \in [0,1]$, by the mean value theorem, there exists $t \in (0,s)$ such that
	\[
		\frac{1}{w^2 - r} \frac{1}{\lvert w^2 - r \rvert^s} = \frac{1}{w^2-r} + s \left[\frac{\partial}{\partial {s'}} \frac{1}{w^2 - r} \frac{1}{\lvert w^2 - r \rvert^{s'}} \right]_{s' = t} =: \frac{1}{w^2-r} + s f(w,r,t).
	\]
	For $r \in [0, D_\gamma/4]$, by applying the mean value theorem to the real and imaginary parts of $f(w,r,t)$, we see that
	\begin{align*}
		\left\lvert \frac{f(w,r,t) - f(w,0,t)}{r} \right\rvert &\ll \max_{0 \leq r \leq D_\gamma/4} \left\lvert \frac{\partial}{\partial r} f(w,r,t) \right\rvert\\
		&\leq \max_{\substack{0 \leq s \leq 1 \\ 0 \leq r \leq D_\gamma/4}} \left\lvert \frac{\partial^2}{\partial r \partial s} \frac{1}{w^2-r} \frac{1}{\lvert w^2 - r \rvert^{s}} \right\rvert \ll \frac{1}{\lvert w \rvert^3}.
	\end{align*}
	Thus, we have
	\[
		\frac{1}{w^2- \frac{D_\gamma}{4a^2}} \frac{1}{\left\lvert w^2 - \frac{D_\gamma}{4a^2} \right\rvert^s} = \frac{1}{w^2-\frac{D_\gamma}{4a^2}} + s \left[\frac{\partial}{\partial s'} \frac{1}{w^2} \frac{1}{\lvert w \rvert^{2s'}} \right]_{s'=t} + O\left(\frac{s}{a^2} \frac{1}{\lvert w \rvert^3} \right).
	\]
	By substituting the equation with the expression~\eqref{Egamma-decomposition}, we divide the sum into three parts.
	
	First, we show that the series
	\[
		y^s \sum_{a \in \mathbb{Z}_{\neq 0}} \frac{1}{\lvert a \rvert^{s+1}} \sum_{\substack{b\ (2\lvert a \rvert)\\ [a,b,\frac{b^2-D_\gamma}{4a}] \sim Q_\gamma}} \sum_{n \in \mathbb{Z}} \frac{1}{\left(z+n+\frac{b}{2a} \right)^2 - \frac{D_\gamma}{4a^2}}
	\]
	converges as $s \to 0$. Since the most inner sum is periodic with period $1$, it has the Fourier expansion of the form
	\[
		\sum_{n \in \mathbb{Z}} \frac{1}{\left(z+n+\frac{b}{2a}\right)^2 - \frac{D_\gamma}{4a^2}} = \sum_{m \in \mathbb{Z}} c_{a,b}(m,y) e^{2\pi i mx},
	\]
	where the Fourier coefficients are given by
	\[
		c_{a,b}(m,y) = e^{\pi imb/a} e^{-2\pi m y} \int_{-\infty+iy}^{\infty+iy}\frac{e^{-2\pi imz}}{z^2 - \frac{D_\gamma}{4a^2}} dz.
	\]
	For $m > 0$, by the inverse Laplace transformation formula~\cite[29.3.57]{AS65} (see also~\cite[(20)]{Koh85}), we have
	\begin{align*}
		\int_{-\infty+iy}^{\infty+iy} \frac{e^{-2\pi imz}}{z^2 - \frac{D_\gamma}{4a^2}} dz = -\frac{2^{3/2} \pi^2 m^{1/2} \lvert a \rvert^{1/2}}{D_\gamma^{1/4}} J_{1/2} \left(\frac{\pi m \sqrt{D_\gamma}}{\lvert a \rvert} \right).
	\end{align*}
	As for $m \leq 0$, since the poles of the integrand are located on the real axis, we can move the integral path as $y \to \infty$ without crossing the poles. By Cauchy's theorem, the integral equals $0$. Therefore, for $\Re(s) > 0$, the first series equals
	\begin{align*}
		-\frac{2^{3/2} \pi^2 }{D_\gamma^{1/4}} y^s \sum_{m=1}^\infty m^{1/2} q^m \sum_{a \in \mathbb{Z}_{\neq 0}} \frac{S_\gamma(a,m)}{\lvert a \rvert^{s+1/2}} J_{1/2} \left(\frac{\pi m \sqrt{D_\gamma}}{\lvert a \rvert} \right).
	\end{align*}
	By \eqref{Hooley-DIT}, the sum uniformly converges on $s \in [0,1/2]$. Therefore, the limit as $s \to 0$ converges to
	\begin{align*}
		&-\frac{2}{\sqrt{D_\gamma}} \sum_{m=1}^\infty \left(\widetilde{\val}_m(\gamma) + 24 \sigma_1(m) \widetilde{\val}_0(\gamma) \right) q^m
	\end{align*}
	by~\cref{val-explicit-formula}.
	
	Second, we show that the series
	\[
		s y^s \sum_{a \in \mathbb{Z}_{\neq 0}} \frac{1}{\lvert a \rvert^{s+1}} \sum_{\substack{b\ (2\lvert a \rvert) \\ [a,b,\frac{b^2-D_\gamma}{4a}] \sim Q_\gamma}} \sum_{n \in \mathbb{Z}} \left[ \frac{\partial}{\partial s'} \frac{1}{\left(z+n+\frac{b}{2a}\right)^2} \frac{1}{\left\lvert z+n+\frac{b}{2a}\right\rvert^{2s'}} \right]_{s'=t}
	\]
	converges as $s \to 0$, where $t \in (0,s)$. We note that interchanging the order of the derivative in $s'$ and the sum over $n \in \mathbb{Z}$ is guaranteed on $s' \in [0,1]$. Again, the Fourier expansion yields that
	\[
		\sum_{n \in \mathbb{Z}} \frac{1}{\left(z+n+\frac{b}{2a}\right)^2} \frac{1}{\left\lvert z+n+\frac{b}{2a}\right\rvert^{2s'}} = \sum_{m \in \mathbb{Z}}  d(m,y,s') e^{\pi imb/a} e^{2\pi imx},
	\]
	where the Fourier coefficients are given by
	\[
		d(m,y,s) = \int_{-\infty}^\infty \frac{e^{-2\pi imx}}{z^2 \left\lvert z\right\rvert^{2s}} dx.
	\]
	By applying the well-known result~\cite[p.353, (5.7)]{Hej83}
	\begin{align*}
		\int_{-\infty}^\infty (1-iY)^{-\alpha} (1+iY)^{-\gamma} e^{iXY} dY = \begin{cases}
			\displaystyle{\frac{2^{\frac{2-\alpha-\gamma}{2}} \pi X^{\frac{\alpha+\gamma}{2} - 1}}{\Gamma(\gamma)} W_{\frac{-\alpha+\gamma}{2}, \frac{\alpha+\gamma-1}{2}} (2X)} &\text{if } X > 0,\\
			\displaystyle{2^{2-\alpha-\gamma} \pi \frac{\Gamma(\alpha+\gamma-1)}{\Gamma(\alpha) \Gamma(\gamma)}} &\text{if } X = 0,
		\end{cases}
	\end{align*}
	for $Y = -\sgn(m)x, X = 2\pi \lvert m \rvert y$ for $m \neq 0$ and $Y = x, X=0$ for $m=0$, we have
	\begin{align*}
		d(m,y,s) &= -y^{-2s-1} \int_{-\infty}^\infty (1-ix)^{-2-s} (1+ix)^{-s} e^{-2\pi imxy} dx\\
		&= \begin{cases}
			- \displaystyle{\frac{y^{-s-1} \pi^{s+1} m^s}{\Gamma(s+2)} W_{1, s+1/2} (4\pi m y)} &\text{if } m>0,\\
			- \displaystyle{y^{-2s-1} 4^{-s} \pi \frac{\Gamma(2s+1)}{\Gamma(s+2) \Gamma(s)}} &\text{if } m=0,\\
			- \displaystyle{\frac{y^{-s-1} \pi^{s+1} \lvert m \rvert^s}{\Gamma(s)} W_{-1, s+1/2} (4\pi \lvert m\rvert y)} &\text{if } m<0,
		\end{cases}
	\end{align*}
	where $W_{\kappa, \mu}(y)$ is the $W$-Whittaker function (see, for instance~\cite[Chapter VII]{MOS66}). Since $d(m,y,s)$ and its derivative in $s$ decays exponentially as $\lvert m\rvert \to \infty$ (see~\cite[Appendix A]{ALR18}), the series we considered is expressed as
	\begin{align}\label{second-series}
		s y^s \sum_{m \in \mathbb{Z}} d_s(m,y,t) e^{2\pi imx} \sum_{a \in \mathbb{Z}_{\neq 0}} \frac{ S_\gamma(a,m)}{\lvert a \rvert^{s+1}} 
	\end{align}
	for $\Re(s) > 0$, where we put $d_s(m,y,s) = \frac{\partial}{\partial s} d(m,y,s)$. We see that the function $d_s(m,y,s)$ is holomorphic at $s=0$. By \eqref{Hooley-DIT} and the definition of the $J$-Bessel function~\eqref{IJ-Bessel}, the above sum over $m \neq 0$ converges absolutely and uniformly on $s \in [0,1]$. Thus the series~\eqref{second-series} without the term for $m=0$ tends to $0$ as $s \to 0$. On the other hand, as for the term for $m=0$,
	\[
		s y^s d_s(0,y,t) \sum_{a \in \mathbb{Z}_{\neq 0}} \frac{S_\gamma(a,0)}{\lvert a \rvert^{s+1}},
	\]
	we see that $d_s(0,y,0) = -\pi/y$ and
	\[
		\lim_{s \to 0} s \sum_{a \in \mathbb{Z}_{\neq 0}} \frac{S_\gamma(a,0)}{\lvert a \rvert^{s+1}} = \frac{6 \widetilde{\val}_0(\gamma)}{\pi^2 \sqrt{D_\gamma}}
	\]
	by \cref{val-explicit-formula}. Therefore, the limit of the second series as $s \to 0$ converges to $-6 \widetilde{\val}_0(\gamma)/\pi y\sqrt{D_\gamma}$.
	
	Third, and finally, the series
	\[
		y^s \sum_{a \in \mathbb{Z}_{\neq 0}} \frac{1}{\lvert a \rvert^{s+1}} \sum_{\substack{b\ (2\lvert a \rvert)\\ [a,b,\frac{b^2-D_\gamma}{4a}] \sim Q_\gamma}} \sum_{n \in \mathbb{Z}} O \left(\frac{s}{a^2} \left\lvert z+n+\frac{b}{2a} \right\rvert^{-3} \right)
	\]
	converges to $0$. In conclusion, we obtain the limit formula
	\[
		\lim_{s \to 0} E_\gamma(z,s) = -\frac{2}{\sqrt{D_\gamma}} \sum_{m=1}^\infty \left(\widetilde{\val}_m(\gamma) + 24 \sigma_1(m) \widetilde{\val}_0(\gamma) \right) q^m - \frac{6\widetilde{\val}_0(\gamma)}{\pi y \sqrt{D_\gamma}},
	\]
	which completes the proof.
\end{proof}

Finally, at the end of this section, we show the following proposition. In fact, this proposition was proved by Duke--Imamo\={g}lu--T\'{o}th~\cite[Theorem 3.3]{DIT17} by using the residue theorem. Here we give an alternative proof. 

\begin{prop}\label{F-modular-gap}
	Let $\gamma \in \Gamma$ be a primitive hyperbolic element with $\sgn(Q_\gamma) > 0$ and $\mathrm{tr}(\gamma) > 2$. For any $\sigma \in \Gamma$, the function $F_\gamma(z)$ defined in \eqref{F-cocycle} satisfies
	\[
		j(\sigma, z)^{-2} F_\gamma(\sigma z) - F_\gamma(z) = \sqrt{D_\gamma} \sum_{\substack{g \in \Gamma_{w_\gamma} \backslash \Gamma \\ w'_{g^{-1} \gamma g} < \sigma^{-1} i\infty < w_{g^{-1} \gamma g}}} \frac{\sgn(Q_{g^{-1} \gamma g})}{Q_{g^{-1} \gamma g}(z,1)}.
	\]
	For $\sigma = T^n$, the right-hand side equals $0$.
\end{prop}

\begin{proof}
	By \cref{limit-formula-for-hyperbolic}, taking the limit of~\cref{Egamma-gap} as $s \to 0$ immediately implies the result, 
	where the terms, including $E_2^*(z)$, vanish since $E_2^*(z)$ is a harmonic Maass form of weight $2$.
\end{proof}


\section{Hyperbolic Rademacher symbol} \label{s4}



\subsection{Definitions} \label{s4-2}


In 2017, Duke--Imamo\={g}lu--T\'{o}th~\cite{DIT17} introduced a hyperbolic analogue of the Rademacher symbol. In this part, we recall their work.

As in \cref{s3-3}, for a primitive hyperbolic element $\gamma= \smat{a & b \\ c & d} \in \Gamma$ with $\sgn(Q_\gamma) := \sgn(c) > 0$ and $\tr(\gamma) > 2$, we consider the generating function of the cycle integrals $F_\gamma(z)$ defined in \eqref{F-cocycle}. As we showed before, $F_\gamma(z)$ is not a modular form but the $i\infty$-component of a locally harmonic (holomorphic) modular form of weight $2$ on $\Gamma$. In particular, \cref{F-modular-gap} asserts that its modular gap gives a weight $2$ rational parabolic cocycle 
\[
	r_\gamma(\sigma, z) := j(\sigma, z)^{-2} F_\gamma(\sigma z) - F_\gamma(z),
\]
whose poles are located at real quadratic irrationals. A function $r: \Gamma \times \bbH \to \mathbb{C}$ is called a weight $k$ \emph{rational cocycle} if $r$ is a rational function on $z$ and $r(\sigma_1 \sigma_2, z) = j(\sigma_2,z)^{-k} r(\sigma_1, \sigma_2 z) + r(\sigma_2, z)$ holds for all $\sigma_1, \sigma_2 \in \Gamma$. Moreover, if $r(T, z) = 0$ identically holds, the cocycle is called a \emph{parabolic cocycle}. Then there uniquely exists the weight $0$ cocycle $R_\gamma(\sigma, z)$ satisfying $\frac{d}{dz} R_\gamma(\sigma, z) = r_\gamma(\sigma, z)$ since both generators $T, S$ of $\Gamma$ are torsion elements. Such cocycle is constructed by $R_\gamma(\sigma, z) := G_\gamma(\sigma z) - G_\gamma(z)$, where $G_\gamma(z)$ is a primitive function
\begin{align}\label{primitive-function}
	G_\gamma(z) := \widetilde{\mathrm{val}}_0(\gamma) z + \frac{1}{2\pi i} \sum_{m=1}^\infty \frac{\widetilde{\mathrm{val}}_m(\gamma)}{m} q^m.
\end{align}
More explicitly, let $L_\gamma(s, \alpha)$ be the Dirichlet series defined by
\[
	L_\gamma(s, \alpha) := \sum_{m=1}^\infty \frac{\widetilde{\val}_m(\gamma) e^{2\pi im\alpha}}{m^s}
\]
with a rational number $\alpha \in \mathbb{Q}$. This series converges for $\Re(s) \gg 0$, and has the meromorphic continuation to $s > 0$ (see~\cite[Theorem 4.1]{DIT17}). Then the cocycle $R_\gamma(\sigma, z)$ for $\sigma = \smat{a & b \\ c & d} \in \Gamma$ with $c \neq 0$ is explicitly given by
\begin{align}\label{Rgamma-cocycle}
	\begin{split}
	R_\gamma(\sigma, z) &= \sum_{\substack{g \in \Gamma_{w_\gamma} \backslash \Gamma \\ w'_{g^{-1} \gamma g} < \sigma^{-1} i\infty < w_{g^{-1} \gamma g}}} \bigg[ \log(z - w_{g^{-1} \gamma g}) - \log(z - w'_{g^{-1} \gamma g}) \bigg]\\
		&\qquad + \frac{1}{2\pi i} L_\gamma(1, a/c) + \widetilde{\val}_0(\gamma) \cdot \frac{a+d}{c},
	\end{split}
\end{align}
where $w_{g^{-1} \gamma g} > w'_{g^{-1} \gamma g}$ are the fixed points of $g^{-1} \gamma g$. Here, we take the principal branch $\Im \log z \in (-\pi, \pi]$. As remarked before, the summation on the right-hand side is a finite sum. For $c = 0$, the cocycle is given by $R_\gamma(\smat{a & b \\ 0 & d},z) = \widetilde{\mathrm{val}}_0(\gamma) b/d$. Under these notations, they defined the hyperbolic analogues of the Dedekind and Rademacher symbol.

\begin{defi}~\cite{DIT17} \label{DIT-def}
	Let $\gamma \in \Gamma$ be a primitive hyperbolic element with $\sgn(Q_\gamma) > 0$ and $\tr(\gamma) > 2$. Then for any $\sigma \in \Gamma$, the hyperbolic Dedekind symbol $\Phi_\gamma(\sigma)$ is defined by
	\[
		\Phi_\gamma(\sigma) := \frac{2}{\pi} \lim_{y \to \infty} \Im R_\gamma(\sigma, iy) = -\frac{1}{\pi^2} \Re L_\gamma(1, \sigma i\infty).
	\]
	If $\sigma \in \Gamma$ is also a hyperbolic element, then the hyperbolic Rademacher symbol $\Psi_\gamma(\sigma)$ is defined by
	\[
		\Psi_\gamma(\sigma) := \lim_{n \to \infty} \frac{\Phi_\gamma(\sigma^n)}{n}.
	\]
\end{defi}


These are analogues of the classical Dedekind and Rademacher symbols in the following sense. For the weight $0$ cocycle $R(\sigma, z) := \log \Delta(\sigma z) - \log \Delta(z)$ given in \eqref{R-definition}, we get a similar characterization of $\Phi(\sigma)$ by
\[
	\Phi(\sigma) = \frac{1}{2\pi} \lim_{y \to \infty} \Im R(\sigma, iy).
\]
Moreover, Rademacher showed in~\cite[Satz 9]{Rad56} for any non-elliptic element $\sigma$ and any integer $n \in \mathbb{Z}$, $\Psi(\sigma^n) = n \Psi(\sigma)$. Thus we also obtain the expression
\[
	\Psi(\sigma) = \lim_{n \to \infty} \frac{\Phi(\sigma^n)}{n}.
\]

The authors~\cite[Theorem 3]{DIT17} established a beautiful connection between $\Psi_\gamma(\sigma)$ and the linking number of two modular knots defined from $\gamma$ and $\sigma$. We here note the geometric interpretations of these symbols they established.

\begin{thm}\cite[Theorem\ 5.2]{DIT17} \label{DIT-Ded}
	For any $\sigma = \smat{a & b \\ c & d} \in \Gamma$ with $c \neq 0$, we have $\Phi_\gamma(-\sigma) = \Phi_\gamma(\sigma^{-1}) = \Phi_\gamma(\sigma)$ and
	\[
		\Phi_\gamma(\sigma) = -\sum_{\substack{g \in \Gamma_{w_\gamma} \backslash \Gamma \\ w'_{g^{-1} \gamma g} < \sigma^{-1} i\infty < w_{g^{-1} \gamma g}}} 1.
	\]
	This counts the number of $\Gamma$-orbits of the geodesic $S_\gamma$ straddling $\sigma^{-1} i\infty \in \mathbb{Q}$, and hence $\Phi_\gamma(\sigma)$ is a non-positive integer. For $\sigma$ with $c = 0$, $\Phi_\gamma(\sigma) = 0$.
\end{thm}

\begin{thm}\cite[Theorem\ 6.4, 6.8]{DIT17} \label{DIT-Rad}
	If both of $\gamma, \sigma$ are primitive hyperbolic elements, and $w_\sigma, w'_\sigma \not\in \{w_{g^{-1} \gamma g} \mid g \in \Gamma_{w_\gamma} \backslash \Gamma\}$, that is, the geodesic $S_{g^{-1} \gamma g}$ never coincide with $S_\sigma$ ignoring their orientations, we have
	\[
		\Psi_\gamma(\sigma) = - \sum_{\substack{g \in \Gamma_{w_\gamma} \backslash \Gamma / \Gamma_{w_\sigma} \\ S_{g^{-1} \gamma g} \text{ intersects with } S_\sigma}} 1.
	\]
	This counts the intersection numbers of two closed geodesics $S_\gamma, S_\sigma$ on $\Gamma \backslash \bbH$, and hence $\Psi_\gamma(\sigma)$ is a non-positive integer.
\end{thm}


\subsection{First explicit formula for $\Psi_\gamma(\sigma)$} \label{s4-3}


Our aims in this section are to prove \cref{main4} and \cref{Main-theorem1}. In \cref{s3}, we studied the Eisenstein series $E_\gamma(z,s)$. \cref{PARABOLIC} asserts that the limit of $E_T(z,s)$ to $s \to 0^+$ gives the harmonic Maass form $E_2^*(z)$. For a hyperbolic element $\sigma \in \Gamma$ with $\sgn(Q_\sigma) > 0$ and $\mathrm{tr}(\sigma) > 2$, its holomorphic part $E_2(z)$ satisfies
\[
	\lim_{n \to \infty} \Re \int_{\sigma^n z_0}^{\sigma^{n+1} z_0} E_2(z) dz = \Psi(\sigma)
\]
as we described in \eqref{Eint-Psi}. In other words, the weight $0$ cocycle $R(\sigma, z) = \log \Delta(\sigma z) - \log \Delta(z)$ satisfies the equation
\[
	\Psi(\sigma) = \frac{1}{2\pi} \lim_{n \to \infty} \Im R(\sigma, \sigma^n z_0)
\]
for any $z_0 \in S_\sigma$. We begin by showing a similar formula for the hyperbolic Rademacher symbol $\Psi_\gamma(\sigma)$.

\begin{thm}\label{Psi-R-red}
	Let $\gamma \in \Gamma$ be a primitive hyperbolic element such that $\sgn(Q_\gamma) > 0, \mathrm{tr}(\gamma) > 2$. Then for a hyperbolic element $\sigma$ and any $z \in \bbH \cup \mathbb{R} \cup \{i\infty\} \setminus \{w_\sigma, w'_\sigma\}$, we have
	\begin{align}\label{Psi-R-red2}
		\Psi_\gamma(\sigma) = \frac{2}{\pi} \lim_{n \to \infty} \Im R_\gamma(\sigma, \sigma^n z).
	\end{align}
\end{thm}

\begin{proof}
	We check the coincidence of both sides by using the explicit formula \eqref{Rgamma-cocycle}. We put $w_\sigma^\infty := \lim_{n \to \infty} \sigma^n z \in \{w_\sigma, w'_\sigma\}$ by \cref{gamma-limit}. Then the right-hand side becomes
	\begin{align*}
		-\frac{1}{\pi^2} \Re L_\gamma(1, \sigma i\infty) +\frac{2}{\pi} \sum_{\substack{g \in \Gamma_{w_\gamma} \backslash \Gamma \\ w'_{g^{-1} \gamma g} < \sigma^{-1} i\infty < w_{g^{-1} \gamma g}}} \left(\arg(w_\sigma^\infty - w_{g^{-1} \gamma g}) - \arg(w_\sigma^\infty - w'_{g^{-1} \gamma g}) \right).
	\end{align*}
	We remark that if $w_\sigma^\infty$ coincides with either of $w_{g^{-1} \gamma g}$ or $w'_{g^{-1} \gamma g}$, that is, $w_{g^{-1} \gamma g} = w_\sigma$ holds, we have to pay careful attention to how the limit is taken. However, such a situation can not happen because $w'_\sigma < \sigma^{-1} i\infty < w_\sigma$ never holds. Therefore, the limit equals
	\begin{align}\label{left-hand}
		\Phi_\gamma(\sigma) + 2 \sum_{\substack{g \in \Gamma_{w_\gamma} \backslash \Gamma \\ w'_{g^{-1} \gamma g} < \sigma^{-1} i\infty, w_\sigma^\infty < w_{g^{-1} \gamma g}}} 1.
	\end{align}
	
	As for the left-hand side, by the definition,
	\begin{align*}
		\Psi_\gamma(\sigma) = \lim_{n \to \infty} \frac{\Phi_\gamma(\sigma^n)}{n} = \lim_{n \to \infty} \frac{1}{n} \cdot \frac{2}{\pi} \lim_{y \to \infty} \Im R_\gamma(\sigma^n, iy).
	\end{align*}
	Since $R_\gamma(\sigma, z)$ satisfies the cocycle condition, that is, $R_\gamma(\sigma_1 \sigma_2, z) = R_\gamma(\sigma_1, \sigma_2 z) + R_\gamma(\sigma_2, z)$, we have
	\[
		R_\gamma(\sigma^n, iy) = \sum_{k=0}^{n-1} R_\gamma(\sigma, \sigma^k iy).
	\]
	Combining them and \eqref{Rgamma-cocycle}, the hyperbolic Rademacher symbol is equal to
	\begin{align*}
		&\frac{2}{\pi} \lim_{n \to \infty} \frac{1}{n} \sum_{k=0}^{n-1} \left(- \frac{1}{2\pi} \Re L_\gamma(1, \sigma i\infty) \right.\\
			&\left. +\sum_{\substack{g \in \Gamma_{w_\gamma} \backslash \Gamma \\ w'_{g^{-1} \gamma g} < \sigma^{-1} i\infty < w_{g^{-1} \gamma g}}} \left(\arg(\sigma^k i\infty - w_{g^{-1} \gamma g}) - \arg(\sigma^k i\infty - w'_{g^{-1} \gamma g}) \right) \right).
	\end{align*}
	For a large enough $k$, the inner sum equals
	\[
		\sum_{\substack{g \in \Gamma_{w_\gamma} \backslash \Gamma \\ w'_{g^{-1} \gamma g} < \sigma^{-1} i\infty, w_\sigma^\infty < w_{g^{-1} \gamma g}}} \pi.
	\]
	Thus $\Psi_\gamma(\sigma)$ coincides with \eqref{left-hand}.
\end{proof}

As corollaries, we get \cref{main4} and \cref{Main-theorem1}.

\begin{cor}[\cref{Main-theorem1}]\label{main-theorem1}
	Let $\gamma, \sigma$ be as in \cref{Psi-R-red}, and $w_\sigma^\infty := \lim_{n \to \infty} \sigma^n z \in \{w_\sigma, w'_\sigma\}$ with $z \in \bbH$. Then we have
	\begin{align*}
		\Psi_\gamma(\sigma) &= \Phi_\gamma(\sigma) + 2 \sum_{\substack{g \in \Gamma_{w_\gamma} \backslash \Gamma \\ w'_{g^{-1} \gamma g} < \sigma^{-1} i\infty, w_\sigma^\infty < w_{g^{-1} \gamma g}}} 1
	\end{align*}
\end{cor}

An important consequence is to make the Rademacher symbol $\Psi_\gamma(\sigma)$ computable. As remarked in~\cite[Remark 3.1]{DIT17}, the sums are finite, and we can easily list all quadratic forms $Q_{g^{-1} \gamma g}$ $(g \in \Gamma_{w_\gamma} \backslash \Gamma)$ such that $w'_{g^{-1} \gamma g} < \sigma^{-1} i\infty < w_{g^{-1} \gamma g}$ by using Mathematica, for example. For convenience, we describe more details. Let $\mathcal{Q}_D$ be the set of binary quadratic forms $[A,B,C]$ with discriminant $D = B^2 - 4AC > 0$. For a given rational number $-d/c \in \mathbb{Q}$ with $c > 0, (c,d) = 1$, we now show the finiteness of the set $\mathcal{Q}_D(-d/c) := \{[A,B,C] \in \mathcal{Q}_D \mid w' < -d/c < w\}$, where $w > w'$ are the roots of $[A,B,C](z,1)$. Since $[A,B,C]$ and $[-A,-B,-C]$ have the same roots and the same discriminant, we always assume that $A > 0$ without loss of generality. The equation
\[
	\frac{-B -\sqrt{D}}{2A} < -\frac{d}{c} < \frac{-B+\sqrt{D}}{2A}
\]
is equivalent to $-c\sqrt{D} < cB - 2dA < c\sqrt{D}$. When we put $n = cB - 2dA \in \mathbb{Z}$, the number of choices of $n$ is finite. Since $4A(d^2A + nd - c^2C) = c^2D-n^2$, the number $c^2 D - n^2$ is a multiple of $4$, and $A$ is a divisor of $(c^2D-n^2)/4$. For each fixed $n$, the number of choices of $A$ is also finite. Now we take $n$ and $A$ as above. If both of
\[
	B = \frac{n + 2dA}{c}, \qquad C = \frac{B^2 - D}{4A}
\]
are integers, the quadratic form $[A,B,C]$ is desired one. This algorithm with finite steps lists all $[A,B,C] \in \mathcal{Q}_D(-d/c)$. By restricting $\mathcal{Q}_D$ to the class $[Q_\gamma] \in \mathcal{Q}_{D_\gamma}/\Gamma$, we can count $\Psi_\gamma(\sigma)$. The classification can be made, for example, by using the continued fraction expansions of fixed points of $Q \in \mathcal{Q}_D(-d/c)$.

\begin{ex}\label{ex1}
	For $\gamma = \sigma = \smat{2 & 1 \\ 1 & 1}$, we compute $\Psi_\gamma(\sigma)$. Note that $\sigma^{-1} i\infty = -1, w_\sigma^\infty = \frac{1+\sqrt{5}}{2}$, and $D_\gamma = 5$. In this case, all quadratic forms in $\mathcal{Q}_5$ are $\Gamma$-equivalent to $Q_\gamma = [1,-1,-1]$. Consider the set $\mathcal{Q}_5^+(-1) := \{[A,B,C] \in \mathcal{Q}_5 \mid A>0, w' < -1 < w\}$. First, we take an integer $n$ such that $-\sqrt{5} < n < \sqrt{5}$ and $5-n^2$ is a multiple of $4$, that is, $n \in \{-1, 1\}$. Second, for each $n$, we take a divisor $A$ of $(5-n^2)/4 = 1$, that is $A = 1$. Thus we have $\mathcal{Q}_5^+(-1) = \{[1,1, -1], [1, 3, 1]\}$. Since $\mathcal{Q}_5(-1) = \mathcal{Q}_5^+(-1) \cup (-\mathcal{Q}_5^+(-1))$, we get $\Phi_\gamma(\sigma) = -4$ by \cref{DIT-Ded}. Moreover, these $4$ quadratic forms never satisfy $w' < w_\sigma^\infty < w$, we also get $\Psi_\gamma(\sigma) = \Phi_\gamma(\sigma) = -4$.
\end{ex}

\begin{ex}\label{ex2}
	For $D = 12$, we know that $\lvert \mathcal{Q}_{12} / \Gamma \rvert = 2$. Let $\sigma = \smat{2 & 1 \\ 1 & 1}$ again, and $\gamma_1 = \smat{3 & 2 \\ 1 & 1}, \gamma_2 = \smat{3& 1\\ 2 & 1}$. In this case, both of $Q_{\gamma_1} = [1,-2,-2]$ and $Q_{\gamma_2} = [2,-2,-1]$ are in $\mathcal{Q}_{12}$, but they are not $\Gamma$-equivalent. By the same argument, we get
	\begin{align*}
		&\{Q \sim Q_{\gamma_1} \mid w' < -1 < w\}\\
		&\quad = \{[-3,-6,-2], [-2,-6,-3], [-2,-2,1], [1,0,-3], [1,2,-2], [1,4,1]\},\\
		&\{Q \sim Q_{\gamma_2} \mid w' < -1 < w\}\\
		&\quad = \{[3,6,2], [2,6,3], [2,2,-1], [-1,0,3], [-1,-2,2], [-1,-4,-1]\},
	\end{align*}
	that is, $\Phi_{\gamma_1}(\sigma) = \Phi_{\gamma_2}(\sigma) = -6$. For the above $12$ quadratic forms, the condition $w' < w_\sigma^\infty < w$ is also satisfied if $Q = [1,0,-3], [-1,0,3]$. Therefore, we get $\Psi_{\gamma_1}(\sigma) = \Psi_{\gamma_2} (\sigma) = -4$.
\end{ex}

\begin{cor}[\cref{main4}] \label{main-theorem4}
	Let $\gamma, \sigma$ be as in \cref{Psi-R-red}. For $z_0 \in \bbH$, we have
	\[
		\Psi_\gamma(\sigma) = 4 \lim_{n \to \infty} \Re \int_{\sigma^n z_0}^{\sigma^{n+1} z_0} \frac{1}{2\pi i} F_\gamma(z) dz.
	\]
\end{cor}

\begin{proof}
	Since $R_\gamma(\sigma, z) = G_\gamma(\sigma z) - G_\gamma(z)$ with $(d/dz) G_\gamma(z) = F_\gamma(z)$ as in \eqref{primitive-function},
	\begin{align*}
		4 \lim_{n \to \infty} \Re \int_{\sigma^n z_0}^{\sigma^{n+1} z_0} \frac{1}{2\pi i} F_\gamma(z) dz = \frac{2}{\pi} \lim_{n \to \infty} \Im R_\gamma(\sigma, \sigma^n z_0).
	\end{align*}
	By \cref{Psi-R-red}, this equals $\Psi_\gamma(\sigma)$.
\end{proof}

The expression of $\Psi_\gamma(\sigma)$ in \cref{DIT-Rad} implies that the hyperbolic Rademacher symbol is a class invariant function if $\sigma$ is primitive with $w_\sigma, w'_\sigma \not\in \{w_{g^{-1} \gamma g} \mid g \in \Gamma_{w_\gamma} \backslash \Gamma\}$. The following proposition removes this assumption.

\begin{prop}\label{Class-inv}
	Under the same assumptions as in \cref{Psi-R-red}, for any $g \in \Gamma$ we have
	\[
		\Psi_\gamma(g^{-1} \sigma g) = \Psi_\gamma(\sigma).
	\]
\end{prop}

\begin{proof}
	By \cref{Psi-R-red},
	\begin{align*}
		\Psi_\gamma(g^{-1} \sigma g) &= \frac{2}{\pi} \lim_{n \to \infty} \Im R_\gamma(g^{-1} \sigma g, g^{-1} \sigma^n g z)\\
			&= \frac{2}{\pi} \lim_{n \to \infty} \Im \bigg[ R_\gamma(g^{-1}, \sigma^{n+1} g z) + R_\gamma(\sigma, \sigma^n g z) + R_\gamma(g, g^{-1} \sigma^n g z)\\
				&\qquad \qquad \qquad \qquad + R_\gamma(g^{-1}, \sigma^n gz) -R_\gamma(g^{-1}, \sigma^n gz) \bigg]. 
	\end{align*}
	The second term becomes $\Psi_\gamma(\sigma)$. The sum of the third and fourth terms is
	\[
		R_\gamma(g, g^{-1} \sigma^n gz) + R_\gamma(g^{-1}, \sigma^ngz) = R_\gamma(g^{-1} g, g^{-1} \sigma^n g z) = 0.
	\]
	The remaining two terms are also cancelled out.
\end{proof}

\begin{rem}
	The expression of $\Psi_\gamma(\sigma)$ in \cref{DIT-Rad} also implies a symmetric property $\Psi_\gamma(\sigma) = \Psi_\sigma(\gamma)$. However, it is not trivial from its definition in \cref{DIT-def}, and the expressions in \cref{Psi-R-red}, \cref{main-theorem1}, and \cref{main-theorem4}. 
\end{rem}


\section{Hyperbolic Dedekind symbol} \label{s5}

\subsection{Cocycle conditions}


It is known that the classical Dedekind symbol $\Phi(\sigma)$ satisfies the following relation~\cite[(62)]{RG72}.
\[
	\Phi(\sigma_1 \sigma_2) = \Phi(\sigma_1) + \Phi(\sigma_2) - 3\sgn(c_1 c_2 c_{12}),
\]
where we put $\sigma_i = \smat{a_i & b_i \\ c_i & d_i} \in \Gamma$ and $\sigma_1\sigma_2 = \smat{a_{12} & b_{12} \\ c_{12} & d_{12}}$. The Dedekind symbol $\Phi(\sigma)$ is uniquely determined by this relation. As for the hyperbolic Dedekind symbol $\Phi_\gamma(\sigma)$, we have the following characterization.

\begin{prop}
	Let $\gamma \in \Gamma$ be a primitive hyperbolic element with $\sgn(Q_\gamma) = \sgn(c) > 0$ and $\mathrm{tr}(\gamma) > 2$. Then for any $\sigma_1, \sigma_2 \in \Gamma$, the hyperbolic Dedekind symbol satisfies
	\[
		\Phi_\gamma(\sigma_1 \sigma_2) = \Phi_\gamma(\sigma_1) + \Phi_\gamma(\sigma_2) + 2 \sum_{\substack{g \in \Gamma_{w_\gamma} \backslash \Gamma \\ w'_{g^{-1} \gamma g} < \sigma_1^{-1} i\infty, \sigma_2 i\infty < w_{g^{-1} \gamma g}}} 1.
	\]
	If either $\sigma_1$ or $\sigma_2$ is in $\Gamma_\infty$, then the last sum equals $0$.
\end{prop}

\begin{proof}
	By the cocycle relation of $R_\gamma(\sigma, z)$,
	\begin{align*}
		\Phi_\gamma(\sigma_1 \sigma_2) &= \frac{2}{\pi} \lim_{y \to \infty} \Im \bigg(R_\gamma(\sigma_1, \sigma_2 iy) + R_\gamma(\sigma_2, iy) \bigg)\\
		&= \frac{2}{\pi} \lim_{y \to \infty} \Im R_\gamma(\sigma_1, \sigma_2 iy) + \Phi_\gamma(\sigma_2).
	\end{align*}
	The explicit formula \eqref{Rgamma-cocycle} implies that
	\begin{align*}
		\frac{2}{\pi} \lim_{y \to \infty} \Im R_\gamma(\sigma_1, \sigma_2 iy) = 2 \sum_{\substack{g \in \Gamma_{w_\gamma} \backslash \Gamma \\ w'_{g^{-1} \gamma g} < \sigma_1^{-1} i\infty, \sigma_2 i\infty < w_{g^{-1} \gamma g}}} 1 - \frac{1}{\pi^2} \Re L_\gamma(1, \sigma_1 i\infty).
	\end{align*}
	Thus we get the relation.
\end{proof}

By this proposition, we can compute the hyperbolic Dedekind symbol $\Phi_\gamma(\sigma)$ inductively. Here we describe fundamental relations for this induction, which have slightly different forms from~\cite[Theorem 5.4]{DIT17}.

\begin{cor}\label{Phi-induction}
	Let $T = \smat{1 & 1 \\ 0 & 1}$ and $S = \smat{0 & -1 \\ 1 & 0}$ be generators of $\Gamma$. For any $\sigma \in \Gamma$, we have
	\begin{align}
		\Phi_\gamma(\sigma^{-1}) &= \Phi_\gamma(-\sigma) = \Phi_\gamma(\sigma), \qquad \Phi_\gamma(I) = 0,\\
		\Phi_\gamma(T^{\pm 1} \sigma) &= \Phi_\gamma(\sigma),\\
		\Phi_\gamma(S^{\pm 1} \sigma) &= \Phi_\gamma(\sigma) + \Phi_\gamma(S) + 2 \sum_{\substack{g \in \Gamma_{w_\gamma} \backslash \Gamma \\ w'_{g^{-1} \gamma g} < 0, \sigma i\infty < w_{g^{-1} \gamma g}}} 1.
	\end{align}
\end{cor}


\subsection{Second explicit formula for $\Psi_\gamma(\sigma)$} \label{s5-2}


In this section, we provide a precise statement of \cref{Main-theorem2} and its proof. Since $\Psi_\gamma(\sigma)$ is a class invariant (\cref{Class-inv}) and satisfies $\Psi_\gamma(\sigma^n) = n \Psi_\gamma(\sigma)$, the problem is reduced to that for primitive elements
\begin{align}\label{gam-sig}
	\gamma = \pmat{a_0 & 1 \\ 1 & 0} \cdots \pmat{a_{2n-1} & 1 \\ 1 & 0}, \quad \sigma = \pmat{b_0 & 1 \\ 1 & 0} \cdots \pmat{b_{2m-1} & 1 \\ 1 & 0}.
\end{align}
Here all coefficients $a_i, b_j$ are positive integers.

First, for a rational number $x \in \mathbb{Q}$, we put
\[
	e_\gamma(x) := \Phi_\gamma(S) + 2 \sum_{\substack{g \in \Gamma_{w_\gamma} \backslash \Gamma \\ w'_{g^{-1} \gamma g} < 0, x < w_{g^{-1} \gamma g}}} 1.
\]
Then \cref{Phi-induction} implies that
\begin{align}\label{TS-rec}
	\Phi_\gamma(T^{\pm 1} \sigma) = \Phi_\gamma(\sigma), \qquad \Phi_\gamma(S^{\pm 1} \sigma) = \Phi_\gamma(\sigma) + e_\gamma(\sigma i\infty).
\end{align}
The idea of our proof is that we first compute the hyperbolic Dedekind symbol $\Phi_\gamma(\sigma)$ explicitly by using \eqref{TS-rec}, and apply the homogenization formula $\Psi_\gamma(\sigma) = \lim_{n \to \infty} \Phi_\gamma(\sigma^n)/n$.

\begin{prop}\label{division}
	For the hyperbolic element $\sigma$ of the form
	\[
		\sigma = \pmat{b_0 & 1 \\ 1 & 0} \cdots \pmat{b_{2m-1} & 1 \\ 1 & 0},
	\]
	we have
	\begin{align*}
		\Phi_\gamma(\sigma) &= \sum_{j =1}^{2m-1} e_\gamma \bigg( (-1)^j [b_j, b_{j+1}, \dots, b_{2m-1}] \bigg) + \Phi_\gamma(S),
	\end{align*}
	where $[b_j, b_{j+1}, \dots, b_{2m-1}] \in \mathbb{Q}$ is a continued fraction expansion $($see \cref{s2}$)$.
\end{prop}

\begin{proof}
	By the fact
	\[
		\pmat{a & 1 \\ 1 & 0} \pmat{b & 1 \\ 1 & 0} = \pmat{1 & a \\ 0 & 1} \pmat{0 & -1 \\ 1 & 0} \pmat{1 & -b \\ 0 & 1} \pmat{0 & 1 \\ -1 & 0} = T^a S T^{-b} S^{-1},
	\]
	the element $\sigma$ is expressed as $\sigma = T^{b_0} S T^{-b_1} S^{-1} \cdots T^{-b_{2m-1}} S^{-1}$. From \eqref{TS-rec}, we see that
	\[
		\Phi_\gamma(\sigma) = \Phi_\gamma(T^{-b_1} S^{-1} T^{b_2} S \cdots T^{-b_{2m-1}} S^{-1}) + e_\gamma(-[b_1, \dots, b_{2m-1}]).
	\]
	From \eqref{TS-rec} again, we have
	\[
		\Phi_\gamma(T^{-b_1} S^{-1} \cdots T^{-b_{2m-1}} S^{-1}) = \Phi_\gamma(T^{b_2} S \cdots T^{-b_{2m-1}} S^{-1}) + e_\gamma([b_2, b_3, \dots, b_{2m-1}]).
	\]
	Repeating this process finishes the proof.
\end{proof}

By virtue of this proposition, our problem is reduced again to getting an explicit formula of $e_\gamma(x)$ for $x = (-1)^j [b_j, b_{j+1}, \dots, b_{2m-1}]$. To this end, it suffices to count the number
\begin{align}\label{count}
	\sum_{\substack{g \in \Gamma_{w_\gamma} \backslash \Gamma \\ w'_{g^{-1} \gamma g} < 0, x < w_{g^{-1} \gamma g}}} 1.
\end{align}
Since $w_\gamma = [\overline{a_0, \dots, a_{2n-1}}]$ is a reduced real quadratic irrational, \cref{orbit} implies
\[
	\Gamma w_\gamma = \bigcup_{r=0}^\infty \bigcup_{\substack{0 \leq i < 2n\\ i \equiv r \ (2)}} A_{r,i},
\]
where
\begin{align*}
	A_{r,i} &= \left\{ [k_0, \dots, k_{r-1}, \overline{a_i, \dots, a_{2n-1}, a_0, \dots, a_{i-1}}]\ \middle\vert\ \begin{array}{l}
		k_0 \in \mathbb{Z}\\
		k_1, \dots, k_{r-1} \in \mathbb{Z}_{>0}\\
		k_{r-1} \neq a_{i-1}
	\end{array}
	\right\}
\end{align*}
Thus the number \eqref{count} is divided as
\begin{align*}
	\sum_{\substack{g \in \Gamma_{w_\gamma} \backslash \Gamma \\ w'_{g^{-1} \gamma g} < 0, x < w_{g^{-1} \gamma g}}} 1 &= \sum_{r = 0}^\infty \sum_{\substack{0 \leq i < 2n \\ i \equiv r \ (2)}} \# \{v \in A_{r, i} \mid v' < 0,x < v \text{ or } v < 0,x < v'\}\\
		&=: \sum_{r=0}^\infty \sum_{\substack{0 \leq i < 2n \\ i \equiv r \ (2)}} N_{r, i} (x),
\end{align*}
where $v'$ is the Galois conjugate of $v$.

\begin{lem}\label{lemma-odd}
	Let $\gamma$ and $\sigma$ be as in \eqref{gam-sig}, and $\ell = 1, \dots, m$. For $x = -[b_{2\ell-1}, \dots, b_{2m-1}]$, we have
	\begin{align*}
		&e_\gamma(x) = -2 \bigg[ \sum_{0 \leq i < 2n} \min(a_i, b_{2\ell-1})\\
			& - \sum_{0 \leq k < n} \bigg( \delta(a_{2k} \geq b_{2\ell-1}) \delta([b_{2\ell}, \dots, b_{2m-1}] \geq [a_{2k-1}, a_{2k-2}, \dots, a_{2k+2\ell-2m}])\\
			&\qquad + \delta(a_{2k-1} \geq b_{2\ell-1}) \delta([b_{2\ell}, \dots, b_{2m-1}] \geq [a_{2k}, a_{2k+1}, \dots, a_{2k-2\ell+2m-1}]) \bigg) \bigg],
	\end{align*}
	where $\delta(P) = 1$ if $P$ is true or empty and $\delta(P) = 0$ otherwise. We put $a_{i} = a_{i \pmod{2n}}$.
\end{lem}

\begin{proof}
	By a simple calculation, for $1 \leq \ell < m$ we have
	\[
		-[b_{2\ell-1}, \dots, b_{2m-1}] = \begin{cases}
			[-b_{2\ell-1}-1, 1, b_{2\ell}-1, b_{2\ell+1}, \dots, b_{2m-1}] &\text{if } b_{2\ell} > 1,\\
			[-b_{2\ell-1}-1, b_{2\ell+1} +1, b_{2\ell+2}, \dots, b_{2m-1}] &\text{if } b_{2\ell} = 1.
		\end{cases}
	\]
	First we count the number of $N_{r,i} = N_{r,i}(x)$. \\
	
	(i) The case of $r = 0$. For $v = [\overline{a_{2k}, \dots, a_{2k+2n-1}}]$, its Galois conjugate is
	\[
		v' = [-1, 1, a_{2k-1}-1, \overline{a_{2k-2}, \dots, a_{2k-2n-1}}]
	\]
	by \cref{CF-conj}. Since $x \leq -1$, the inequation $v' < x < 0 < v$ never holds. Thus $N_{0, 2k} = 0$ for any $0 \leq k < n$.
	
	(ii) The case of $r=1$. For $v = [k_0, \overline{a_{2k+1}, \dots, a_{2k+2n}}]$ with $k_0 \neq a_{2k}$, we have
	\[
		v' = [k_0-a_{2k}-1, 1, a_{2k-1} - 1, \overline{a_{2k-2}, \dots, a_{2k-2n-1}}].
	\]
	Now we consider the inequation $v' < x < 0 < v$. From $v > 0$, the number $k_0$ should be $k_0 \geq 0$. The inequation $v' < x$ holds if and only if $k_0 - a_{2k}-1 < -b_{2\ell-1}-1$, or $k_0 = a_{2k}-b_{2\ell-1}$ and $[b_{2\ell}, \dots, b_{2m-1}] \geq [a_{2k-1}, a_{2k-2}, \dots, a_{2k+2\ell-2m}]$. Thus the number of such $k_0$ is given by
	\begin{align}\label{odd-1}
		\begin{split}
		N_{1, 2k+1} &= a_{2k} - \min(a_{2k}, b_{2\ell-1})\\
			&\quad + \delta(a_{2k} \geq b_{2\ell-1}) \delta([b_{2\ell}, \dots, b_{2m-1}] \geq [a_{2k-1}, a_{2k-2}, \dots, a_{2k+2\ell-2m}]).
		\end{split}
	\end{align}
	If $\ell = m$, then $x = -b_{2m-1}$. In this case, $k_0$ should satisfy $0 \leq k_0 < a_{2k} - b_{2m-1} + 1$. Thus $N_{1, 2k+1} = a_{2k} - \min(a_{2k}, b_{2m-1}) + \delta(a_{2k} \geq b_{2m-1})$, which is contained in \eqref{odd-1}.
	
	(iii) The case of $r=2$. For $v = [k_0, k_1, \overline{a_{2k+2}, \dots, a_{2k+2n+1}}]$ with $k_1 \neq a_{2k+1}$, its Galois conjugate is given by
	\[
		v' = \left\{\begin{array}{ll}
			[k_0, \dots] &\text{if } k_1 > a_{2k+1} + 1,\\
			\left[k_0+1, \dots \right] &\text{if } k_1 = a_{2k+1} + 1, a_{2k} > 1,\\
			\left[k_0+a_{2k-1}+1, \overline{a_{2k-2}, \dots, a_{2k-2n-1}}\right] &\text{if } k_1 = a_{2k+1} + 1, a_{2k} =1,\\
			\left[ k_0 -1, \dots \right] &\text{if } 0 < k_1 < a_{2k+1}.
		\end{array} \right.
	\]
	Since $\lvert x - 0\rvert \geq 1$, the contribution comes only from the third case. Then we consider $v < x < 0 < v'$. From $v' > 0$, we have $k_0 \geq -a_{2k-1}-1$. On the other hand, from $v < x$, we get $k_0 < -b_{2\ell-1}-1$, or $k_0 = -b_{2\ell-1}-1$ and $[a_{2k+1}+1, a_{2k+2}, \dots, a_{2k-2\ell+2m+1}] \geq [1, b_{2\ell}-1, b_{2\ell+1}, \dots, b_{2m-1}]$. If $b_{2\ell} > 1$, then the last inequality always holds, but for $b_{2\ell} = 1$, this is equivalent to $[a_{2k+1}, \dots, a_{2k-2\ell+2m-1}] \geq [b_{2\ell+1}, \dots, b_{2m-1}]$. Thus we have
	\begin{align}\label{odd-2}
	\begin{split}
		N_{2, 2k+2} 
			&= \delta(a_{2k}=1) \bigg(a_{2k-1} - \min(a_{2k-1}, b_{2\ell-1}) \\
			&\quad + \delta(a_{2k-1} \geq b_{2\ell-1}) \delta([b_{2\ell}, \dots, b_{2m-1}] \geq [a_{2k}, \dots, a_{2k-2\ell+2m-1}]) \bigg).
		\end{split}
	\end{align}
	The case of $\ell = m$ is also included.
	
	(iv) The case of $r = 3$. For $v = [k_0, k_1, k_2, \overline{a_{2k+1}, \dots, a_{2k+2n}}]$ with $k_2 \neq a_{2k}$, its Galois conjugate is given by
	\[
		v' = \left\{\begin{array}{ll}
			[k_0, \dots] &\text{if } k_2 > a_{2k},\\
			\left[k_0, \dots \right] &\text{if } 0 < k_2 < a_{2k}, k_1 > 1,\\
			\left[ k_0 + a_{2k-1} + 1, \overline{a_{2k-2}, \dots, a_{2k-2n-1}}\right] &\text{if } 0 < k_2 = a_{2k} -1, k_1 = 1,\\
			\left[ k_0+1, \dots \right] &\text{if } 0 < k_2 < a_{2k} -1, k_1 = 1.
		\end{array} \right.
	\]
	In this case also, the contribution comes only from the third case. Then we consider $v < x < 0 < v'$. From $v' > 0$, the number $k_0$ should be $k_0 \geq -a_{2k-1} -1$. The inequation $v < x$ is equivalent to $k_0 < -b_{2\ell-1}-1$, or $k_0 = -b_{2\ell-1}-1$ and $[b_{2\ell}, \dots, b_{2m-1}] \geq [a_{2k}, a_{2k+1}, \dots, a_{2k-2\ell+2m-1}]$. Thus we have
	\begin{align}\label{odd-3}
		\begin{split}
		N_{3, 2k+1} &= \delta(a_{2k} > 1) \bigg( a_{2k-1} - \min(a_{2k-1}, b_{2\ell-1})\\
			&\quad + \delta(a_{2k-1} \geq b_{2\ell-1}) \delta([b_{2\ell}, \dots, b_{2m-1}] \geq [a_{2k}, a_{2k+1}, \dots, a_{2k-2\ell+2m-1}]) \bigg).
		\end{split}
	\end{align}
	
	For $r \geq 4$, we easily see that the number $N_{r, i}$ always becomes $0$. Combining all results, we get
	\begin{align*}
		&\sum_{\substack{g \in \Gamma_{w_\gamma} \backslash \Gamma \\ w'_{g^{-1} \gamma g} < 0, x < w_{g^{-1} \gamma g}}} 1 = \sum_{k=1}^n \bigg( a_{2k-1} + a_{2k} - \min(a_{2k-1}, b_{2\ell-1}) -\min(a_{2k}, b_{2\ell-1})\\
			&\qquad + \delta(a_{2k-1} \geq b_{2\ell-1}) \delta([b_{2\ell}, \dots, b_{2m-1}] \geq [a_{2k}, a_{2k+1}, \dots, a_{2k-2\ell+2m-1}])\\
			&\qquad + \delta(a_{2k} \geq b_{2\ell-1}) \delta([b_{2\ell}, \dots, b_{2m-1}] \geq [a_{2k-1}, a_{2k-2}, \dots, a_{2k+2\ell-2m}]) \bigg).
	\end{align*}
	
	To finish the proof, we next compute
	\[
		\Phi_\gamma(S) = - \sum_{\substack{g \in \Gamma_{w_\gamma} \backslash \Gamma \\ w'_{g^{-1} \gamma g} < 0 < w_{g^{-1} \gamma g}}} 1
	\]
	in a similar way. By repeating the same argument as above for $x = 0$, we obtain (i) $N_{0, 2k} = 1$, (ii) $N_{1, 2k+1} = a_{2k}$, (iii) $N_{2, 2k+2} = \delta(a_{2k} > 1) + \delta(a_{2k} = 1) (a_{2k-1} + 1) + a_{2k+1}-1$, (iv) $N_{3, 2k+1} = \delta(a_{2k} > 1) (a_{2k-1} + 1) + a_{2k}-2+\delta(a_{2k} = 1)$, and $N_{r, i} = 0$ for $r \geq 4$. Therefore, we get
	\begin{align}\label{Phi-simple}
		\Phi_\gamma(S) = -2 \sum_{k=1}^n \bigg(a_{2k-1} + a_{2k} \bigg).
	\end{align}
	This concludes the proof.
\end{proof}

\begin{cor}
	For a hyperbolic element $\gamma$ of the form
	\[
		\gamma = \pmat{a_0 & 1 \\ 1 & 0} \cdots \pmat{a_{2n-1} & 1 \\ 1 & 0},
	\]
	we have
	\[
		\Phi_\gamma(S) = -2 \sum_{0 \leq i < 2n} a_i.
	\]
\end{cor}

The corollary is an interesting analogue of the classical Rademacher symbol in the sense of the expression~\eqref{Psi-explicit}. In particular, this gives the number of simple quadratic forms in the $\Gamma$-equivalent class of $Q_\gamma$. Here $Q(X,Y) = AX^2+ BXY + CY^2$ is said to be \emph{simple} if $AC < 0$ holds, which was studied by Choie--Zagier~\cite{CZ93}.

Similarly, we get a formula of $e_\gamma(x)$ for $x = [b_{2\ell}, b_{2\ell+1}, \dots, b_{2m-1}]$ as follows.

\begin{lem}\label{lemma-even}
	Let $\gamma$ and $\sigma$ be as in \eqref{gam-sig}, and $\ell = 1, \dots, m-1$. For $x = [b_{2\ell}, \dots, b_{2m-1}]$, we have
	\begin{align*}
		&e_\gamma(x) = -2 \bigg[ \sum_{0 \leq i < 2n} \min(a_i, b_{2\ell})\\
			& - \sum_{0 \leq k < n} \bigg( \delta(a_{2k} \geq b_{2\ell}) \delta([b_{2\ell+1}, \dots, b_{2m-1}] > [a_{2k+1}, a_{2k+2}, \dots, a_{2k-2\ell+2m-1}])\\
			&\qquad + \delta(a_{2k-1} \geq b_{2\ell}) \delta([b_{2\ell+1}, \dots, b_{2m-1}] > [a_{2k-2}, a_{2k-3}, \dots, a_{2k+2\ell-2m}]) \bigg) \bigg].
	\end{align*}
\end{lem}

Combining \cref{division} with \cref{lemma-odd} and \cref{lemma-even}, we get an explicit formula for $\Phi_\gamma(\sigma)$. Furthermore, by using the homogenization formula
\[
	\Psi_\gamma(\sigma) = \lim_{n \to \infty} \frac{\Phi_\gamma(\sigma^n)}{n},
\]
we get the following formula for the hyperbolic Rademacher symbol $\Psi_\gamma(\sigma)$. 

\begin{thm}[\cref{Main-theorem2}] \label{explicit-psi}
	Let $\gamma, \sigma$ be as in \eqref{gam-sig}. Then we have
	\begin{align*}
		\Psi_\gamma(\sigma) = -2 \left( \sum_{\substack{0 \leq i < 2n \\ 0 \leq j < 2m}} \min(a_i, b_j) - \psi_\gamma(\sigma) \right),
	\end{align*}
	where $\psi_\gamma(\sigma)$ is given by
	\begin{align*}
		\psi_\gamma(\sigma) &= \sum_{0 \leq k < n} \sum_{0 \leq \ell < m} \bigg( \delta(a_{2k} \geq b_{2\ell-1}) \delta([\overline{b_{2\ell}, \dots, b_{2\ell+2m-1}}] \geq [\overline{a_{2k-1}, \dots, a_{2k-2n}}])\\
			&\qquad + \delta(a_{2k-1} \geq b_{2\ell-1}) \delta([\overline{b_{2\ell}, \dots, b_{2\ell+2m-1}}] \geq [\overline{a_{2k}, \dots, a_{2k+2n-1}}])\\
			&\qquad + \delta(a_{2k} \geq b_{2\ell}) \delta([\overline{b_{2\ell+1}, \dots, b_{2\ell+2m}}] > [\overline{a_{2k+1} ,\dots, a_{2k+2n}}])\\
			&\qquad + \delta(a_{2k-1} \geq b_{2\ell}) \delta([\overline{b_{2\ell+1}, \dots, b_{2\ell+2m}}] > [\overline{a_{2k-2}, \dots, a_{2k-2n-1}}]) \bigg).
	\end{align*}
	Here we put $a_i = a_{i \pmod{2n}}$ and $b_j = b_{j \pmod{2m}}$.
\end{thm}

\begin{proof}
	We can also apply \cref{lemma-odd} and \cref{lemma-even} to 
	\[
		\sigma^r = \pmat{b_0 & 1 \\ 1 & 0} \cdots \pmat{b_{2mr-1} & 1 \\ 1 & 0}
	\]
	with $b_j = b_{j \pmod{2m}}$. We can see that the value $\delta([b_{2\ell}, \dots, b_{2mr-1}] \geq [a_{2k-1}, a_{2k-2}, \dots, a_{2k+2\ell-2mr}])$ is equal to
	\begin{align}\label{eq-period}
		\delta([\overline{b_{2\ell}, \dots, b_{2\ell + 2m-1}}] \geq [\overline{a_{2k-1}, \dots, a_{2k-2n}}])
	\end{align}
	if the length $2mr - 2\ell \geq \mathrm{lcm}(2n, 2m)$. The point is that $\mathrm{lcm}(2n,2m)$ and \eqref{eq-period} are independent of $r$. Similar as above, we can estimate other three terms. Thus the limit
	\[
		\Psi_\gamma(\sigma) = \lim_{r \to \infty} \bigg[ \frac{1}{r} \sum_{j=1}^{2mr-1} e_\gamma \bigg((-1)^j [b_j, b_{j+1}, \dots, b_{2mr-1}] \bigg) + \frac{\Phi_\gamma(S)}{r} \bigg]
	\]
	gives our theorem.
\end{proof}

\begin{rem}
	The integer-valued function $\psi_\gamma(\sigma)$ satisfies $0 \leq \psi_\gamma(\sigma) \leq 2mn$, whose bound $2mn$ is independent of the coefficients of $\gamma$ and $\sigma$. We do not know if there exists a simpler formula for the $\psi_\gamma(\sigma)$.
\end{rem}

\begin{ex}
	We compute $\Psi_\gamma(\sigma)$ for $\gamma = \smat{3 & 2 \\ 1 & 1}$ and $\sigma = \smat{2 & 1 \\ 1& 1}$ by using \cref{explicit-psi}. As we showed in \cref{ex2}, we have $\Psi_\gamma(\sigma) = -4$. First, we express these matrices as
	\[
		\gamma = \pmat{2 & 1 \\ 1 & 0} \pmat{1& 1 \\ 1 & 0}, \qquad \sigma = \pmat{1 & 1 \\ 1 & 0} \pmat{1 & 1 \\ 1 & 0}.
	\]
	Then we have $\Psi_\gamma(\sigma) = -2(4 - \psi_\gamma(\sigma))$ with
	\begin{align*}
		\psi_\gamma(\sigma) &= \delta(2 \geq 1) \delta([\overline{1,1}] \geq [\overline{1,2}]) + \delta(1 \geq 1) \delta([\overline{1,1}] \geq [\overline{2,1}])\\
			&\quad + \delta(2 \geq 1) \delta([\overline{1,1}] > [\overline{1,2}]) + \delta(1 \geq 1) \delta([\overline{1,1}] > [\overline{2,1}])\\
			&= 2.
	\end{align*}
	Here note that $[\overline{1,1}] = (1+\sqrt{5})/2, [\overline{1,2}] = (1+\sqrt{3})/2$, and $[\overline{2,1}] = 1+\sqrt{3}$. Thus we get $\Psi_\gamma(\sigma) = -4$ again.
\end{ex}

\appendix
\def\thesection{\Alph{section}}

\section{Continued fractions} \label{s2} 

In this Appendix, we review the theory of continued fractions. It is well-known that a real number $w \in \mathbb{R} \setminus \mathbb{Q}$ has the following unique expression
\[
	w = [a_0, a_1, a_2, \dots] = a_0 + \cfrac{1}{a_1 + \cfrac{1}{a_2 + \cfrac{1}{\ddots}}}
\]
where $a_j \in \mathbb{Z}$ with $a_j \geq 1$ for $j \geq 1$. This expression is called the (simple) \emph{continued fraction expansion} of $w$. An ancient result of Euler and Lagrange asserts that a continued fraction of $w \in \mathbb{R} \backslash \mathbb{Q}$ is periodic if and only if it is a real quadratic irrationality. For simplicity, we write $w = [k_0, \dots, k_{r-1}, \overline{a_0, \dots, a_{n-1}}] := [k_0, \dots, k_{r-1}, a_0, \dots, a_{n-1}, a_0, \dots, a_{n-1}, \dots]$. In addition, this expansion is purely periodic precisely if and only if $w$ is \emph{reduced}, that is, $w>1, -1 < w' < 0$ for the Galois conjugate $w'$ of $w$.

For any real quadratic irrational number $w = [k_0, \dots, k_{r-1}, \overline{a_0, \dots, a_{n-1}}]$, we define $\delta_w \in \SL_2(\mathbb{Z})$ by
\[
	\delta_w := \left\{\begin{array}{ll}
		\pmat{k_0&1\\1&0} \cdots \pmat{k_{r-1}&1\\1&0} &\text{if $r$ is even},\\
		\ \\
		\pmat{k_0&1\\1&0} \cdots \pmat{k_{r-1}&1\\1&0} \pmat{a_0&1\\1&0} \quad &\text{if $r$ is odd}.	\end{array} \right.
\]
Then $\delta_w^{-1} w$ is also a real quadratic irrational, and has a purely periodic expansion. For compatibility with the actions of $\SL_2(\mathbb{Z})$ instead of $\GL_2(\mathbb{Z})$, we can assume that $2r$ and $2n$ are minimal even integers such that $w = [k_0, \dots, k_{2r-1}, \overline{a_0, \dots, a_{2n-1}}]$ without loss of generality. For instance, $w = (3+\sqrt{5})/2$ has the continued fraction expansion $w = [2,1, \overline{1,1}]$, not $w = [2, \overline{1}]$ here.

Let $w = [k_0, \dots, k_{2r-1}, \overline{a_0, \dots, a_{2n-1}}]$ be a real quadratic irrational with minimal $r, n \in \mathbb{Z}$, and $\Gamma = \SL_2(\mathbb{Z})$. We now consider the stabilizer subgroup
\[
	\Gamma_w := \left\{ \gamma = \pmat{a & b \\ c & d} \in \Gamma \mid \gamma w = \frac{aw+b}{cw+d} = w \right\}.
\] 
If we take $\delta_w \in \Gamma$ as above, then $\delta_w^{-1} \Gamma_w \delta_w = \Gamma_{\delta_w^{-1} w}$ holds. For the reduced $\delta_w^{-1} w$, it is well-known that
\[
	\Gamma_{\delta_w^{-1} w} := \left\{ \pm \gamma_w^n \mid n \in \mathbb{Z}\right\},
\]
where
\[
	\gamma_w := \pmat{a_0 & 1 \\ 1 & 0} \cdots \pmat{a_{2n-1} & 1 \\ 1 & 0} \in \Gamma
\]
is defined from the period of $w$. Thus we reach the following lemma.

\begin{lem}\label{stab}
	For a real quadratic irrational $w = [k_0, \dots, k_{2r-1}, \overline{a_0, \dots, a_{2n-1}}]$ with minimal $r,n \in \mathbb{Z}$, we put
	\[
		\delta_w := \pmat{k_0&1\\1&0} \cdots \pmat{k_{2r-1}&1\\1&0}, \quad \gamma_w := \pmat{a_0 & 1 \\ 1 & 0} \cdots \pmat{a_{2n-1} & 1 \\ 1 & 0} \in \Gamma.
	\]
	Then the stabilizer $\Gamma_w$ is given by
	\[
		\Gamma_w = \left\{ \pm \delta_w \gamma_w^n \delta_w^{-1} \mid n \in \mathbb{Z} \right\}.
	\]
\end{lem}

Moreover, from the easy facts that
\[
	\pmat{k & 1 \\ 1 & 0} [a_0, a_1, \dots] = [k, a_0, a_1, \dots] \quad \text{and} \quad \det \pmat{k & 1 \\ 1 & 0} = -1,
\]
we have the following lemma.

\begin{lem}\label{orbit}
	For a reduced real quadratic irrational number $w = [\overline{a_0, \dots, a_{2n-1}}]$ with a minimal even length $2n \in 2\mathbb{Z}$, the $\Gamma$-orbit of $w$ is given by
	\[
		\Gamma w = \bigcup_{r=0}^\infty \bigcup_{\substack{0 \leq i < 2n\\ i \equiv r \ (2)}} A_{r,i},
	\]
	where
	\begin{align*}
		A_{r,i} &= \{ [k_0, \dots, k_{r-1}, \overline{a_i, \dots, a_{2n-1}, a_0, \dots, a_{i-1}}] \mid\\
			&\qquad k_0 \in \mathbb{Z}, k_1, \dots, k_{r-1} \in \mathbb{Z}_{>0}, k_{r-1} \neq a_{i-1}\}.
	\end{align*}
\end{lem}

Finally, we show a useful proposition on the continued fraction expansions of the Galois conjugate $w'$. To give it, we prepare an easy lemma.

\begin{lem}\label{CF-lem}
	Let $a, k \in \mathbb{Z}, b \in \mathbb{Z}_{\neq -1}$ be integers, and $w \in \mathbb{R}$. Then we have
	\begin{align*}
		\pmat{0&-1\\1&0} w &= [-1, 1, w-1],\\
		\pmat{k&1\\1&0} [-1, 1, a-1, w] &= [k-a-1, 1, w-1],\\
		\pmat{k&1\\1&0} [b, 1, w-1] &= \left\{\begin{array}{ll}
			[k, b, 1, w-1] &\text{if } b >-1,\\ 
			\left[ k-1, 1, -b-2, w \right] &\text{if } b < -1.
		\end{array} \right.
	\end{align*}
\end{lem}

\begin{proof}
	We here give a proof for the first formula. By the definition, we see that
	\[
		\pmat{0&-1\\1&0} w = -\frac{1}{w} = -1 + \frac{w-1}{w} = -1 + \cfrac{1}{1 + \cfrac{1}{w-1}} = [-1, 1, w-1].
	\]
	Similarly, we can obtain the remaining formulas.
\end{proof}

By using this lemma, we can compute the continued fraction expansion of the Galois conjugate $w'$ for $w = [k_0, \dots, k_{r-1}, \overline{a_0, \dots, a_{2n-1}}]$ with $k_{r-1} \neq a_{2n-1}$. By the well-known fact~\cite[Lemma 1.28]{Aig13} that
\[
	S\alpha' = -\frac{1}{\alpha'} = [\overline{a_{t-1}, \dots, a_1, a_0}]
\]
for a reduced real quadratic irrational $\alpha = [\overline{a_0, a_1, \dots, a_{t-1}}]$, we have
\[
	w' = \pmat{k_0 & 1 \\ 1 & 0} \cdots \pmat{k_{r-1} & 1 \\ 1 & 0} S^{-1} \cdot [\overline{a_{2n-1}, \dots, a_0}].
\]
When $r = 0$, by the first formula of \cref{CF-lem},
\[
	w' = S^{-1} \cdot [\overline{a_{2n-1}, \dots, a_0}] = [-1, 1, a_{2n-1}-1, \overline{a_{2n-2}, \dots, a_0, a_{2n-1}}].
\]
When $r = 1$, by the second formula of \cref{CF-lem},
\[
	w' = [k_0 - a_{2n-1} -1, 1, a_{2n-2}-1, \overline{a_{2n-3}, \dots, a_0, a_{2n-1}, a_{2n-2}}].
\]
When $r=2$, by the third formula of \cref{CF-lem},
\[
	w' = [k_0, k_1 - a_{2n-1} -1, 1, a_{2n-2}-1, \overline{a_{2n-3}, \dots, a_0, a_{2n-1}, a_{2n-2}}]
\]
if $k_1 > a_{2n-1}$, and 
\[
	w' = [k_0-1, 1, a_{2n-1} -k_1 -1, \overline{a_{2n-2}, \dots, a_0, a_{2n-1}}]
\]
if $0 < k_1 < a_{2n-1}$. Repeating the process for $r \geq 3$, we obtain the proposition.

\begin{prop}\label{CF-conj}
	Let $w = [k_0, \dots, k_{r-1}, \overline{a_0, \dots, a_{2n-1}}]$ be a real quadratic irrational with $k_{r-1} \neq a_{2n-1}$. For any integer $j \in \mathbb{Z}$, we put $a_j = a_{j \pmod{2n}}$. Then its Galois conjugate $w'$ has the following continued fraction expansion.
	\begin{align*}
		w' = \left\{\begin{array}{ll}
			[-1, 1, a_{2n-1}-1, \overline{a_{2n-2}, \dots, a_{-1}}] &\text{if } r = 0,\\
			\left[k_0 - a_{2n-1} -1, 1, a_{2n-2}-1, \overline{a_{2n-3}, \dots, a_{-2}}\right] &\text{if } r = 1,\\
			\left[k_0, \dots, k_{r-2}, k_{r-1} - a_{2n-1} -1, 1, a_{2n-2}-1, \overline{a_{2n-3}, \dots, a_{-2}}\right] &\text{if } r \geq 2,\\
			\left[k_0, \dots, k_{r-3}, k_{r-2} -1, 1, a_{2n-1} -k_{r-1} -1, \overline{a_{2n-2}, \dots, a_{-1}} \right] &\text{if } r \geq 2,
		\end{array} \right.
	\end{align*}
	where the third line holds if $k_{r-1} > a_{2n-1}$, and the fourth line holds if $0 < k_{r-1} < a_{2n-1}$. Moreover, if an inner entry becomes $0$, then we regard as
	\[
		[\dots, a, 0, b, \dots] = [\dots, a+b, \dots].
	\]
\end{prop}

\begin{ex}
	Let $w = [2, 1, \overline{1,4,3,2}] = \frac{36 + 2\sqrt{39}}{19}$. Then by the fourth line of \cref{CF-conj}, we have
	\[
		w' = [2 - 1, 1, 2-1-1, \overline{3,4,1,2}] = [1, 1, 0, \overline{3,4,1,2}].
	\]
	Now a $0$-entry appears so that this equals
	\[
		w' = [1, 4, \overline{4,1,2,3}] = \frac{36 - 2\sqrt{39}}{19}.
	\]
\end{ex}

\section{Hooley's estimate}\label{Appendix-B}

Suppose that $\gamma = \smat{a & b \\ c & d} \in \Gamma$ is a primitive hyperbolic element with $a+d > 2$. Let $S_\gamma(k,m)$ be the exponential sum defined in \eqref{exponential-sum}. The aim of this section is to recall and complete the proof of following Hooley's estimate, which was omitted in~\cite{Hoo63}.

\begin{thm}[{Hooley~\cite{Hoo63}}]\label{Hooley-theorem1}
	For $m > 0$ and $x > 0$, we have for any $\epsilon > 0$ that
	\[
		\sum_{0 < k \leq x} S_\gamma(k,m) \ll_{\gamma, \epsilon} m^\epsilon x^\epsilon(m^{4/5} + x^{3/4}).
	\]
\end{thm}

By \cref{stab} and an easy fact that $S_{g^{-1} \gamma g}(k,m) = S_\gamma(k,m)$ for any $g \in \Gamma$, we can assume that the primitive hyperbolic element $\gamma$ is either of the form
\[
	\gamma = \pmat{a& b \\ c & d} = \pmat{a_0 & 1 \\ 1 & 0} \cdots \pmat{a_{2n-1} & 1 \\ 1 & 0}
\]
or its inverse matrix. We prove here \cref{Hooley-theorem1} only for the former case, but the latter case can be proved as well. In this case, we see that $a, b, c$, and $d$ are all positive, and $a \geq c > 0$ holds.

Let $k > 0$ be a positive integer. For any pair $(r, s) \in \mathbb{Z}^2$ with $Q_\gamma(r,s) = k$ and $\gcd(r,s) = 1$, the equation $r\sigma - s\rho = 1$ has the solution
\[
	(\rho, \sigma) = (\rho_t, \sigma_t) = (\rho_0 + rt, \sigma_0 + st) \quad (t \in \mathbb{Z}).
\]
For the solution, we put
\[
	g_t = \pmat{r & \rho_t \\ s & \sigma_t} = \pmat{r & \rho_0 \\ s & \sigma_0} T^t \in \Gamma,
\]
where $T = \smat{1 & 1 \\ 0 & 1}$. Then we have
\begin{align*}
	(Q_\gamma \circ g_t)(X,Y) &= [Q_\gamma(r,s), 2cr\rho_t + (d-a)(r\sigma_t + s\rho_t) -2b s\sigma_t, *]\\
		&=: [k, \ell_t, \frac{\ell_t^2 - D_\gamma}{4k}].
\end{align*}
Since $\ell_t = \ell_0 + 2Q_\gamma(r,s) t$ holds, there uniquely exists $t \in \mathbb{Z}$ such that $0 \leq \ell_t < 2k$. The correspondence $(r,s) \mapsto \ell_t$ induces the bijection
\[
	\left\{(r,s) \in \mathbb{Z}^2\ \middle\vert\ \begin{array}{l}
		Q_\gamma(r,s) = k\\
		r > 0, 0 < s \leq \frac{c}{a}r \\ 
		(r,s) = 1
		\end{array}
		\right\} \to \left\{\ell \in [0, 2k)\ \middle\vert\ [k, \ell, \frac{\ell^2 - D_\gamma}{4k}] \sim Q_\gamma \right\},
\]
where $\{(X,Y) \in \mathbb{R}^2 \mid X>0, 0 < Y \leq (c/a)X\}$ is a fundamental domain of $\Gamma_{w_\gamma} = \{\pm \gamma^n \mid n \in \mathbb{Z}\}$ in the region $\{(X,Y) \in \mathbb{R}^2 \mid Q_\gamma(X,Y) > 0\}$. Under the correspondence, we have
\begin{align*}
	\ell &\equiv \frac{-2 \overline{s} k +(d-a)r - 2bs}{r} \pmod{2k}
\end{align*}
with $s \overline{s} \equiv 1 \pmod{r}$. Therefore, we see that
\begin{align*}
	\sum_{0 < k \leq x} S_\gamma(k,m) = \sum_{0 < k \leq x} \sum_{\substack{0 \leq \ell < 2k \\ [k,\ell,\frac{\ell^2-D_\gamma}{4k}] \sim Q_\gamma}} e^{\pi i m \ell/k} = \sum_{\substack{r > 0, 0 < s \leq \frac{c}{a} r\\ 0 < Q_\gamma(r,s) \leq x\\ (r,s) = 1}} e^{-2\pi i m \overline{s}/r} \varphi_m(r,s),
\end{align*}
where we put
\[
	\varphi_m(r,s) = \exp \left(2\pi i m \frac{(d-a)r - 2bs}{2r(cr^2 + (d-a)rs - bs^2)} \right).
\]
The expression is also mentioned in~\cite[p.110]{Hoo63}. By the conditions assumed to $(r,s)$, we can show that there exists constants $A, B > 0$ depending on $\gamma$ such that $r \leq A x^{1/2}$ and $cr^2 + (d-a)rs - bs^2 \geq Br^2$. Then the above sum is also expressed as
\begin{align}\label{Sgamma-expre}
	\sum_{0 < k \leq x} S_\gamma(k,m) = \sum_{0 < r \leq A x^{1/2}} \sum_{\substack{0 < s \leq \frac{c}{a} r \\ 0 < Q_\gamma(r,s) \leq x \\ (r,s) = 1}} e^{-2\pi i m\overline{s}/r} \varphi_m(r,s).
\end{align}
Here we note that the set $I = \{s \in \mathbb{R} \mid 0 < s \leq (c/a)r, 0 < Q_\gamma(r,s) \leq x\}$ is an interval or a union of two intervals contained in $(0, cr/a] \subset [0,r]$ by $a \geq c > 0$.

We now evaluate the sum~\eqref{Sgamma-expre}. First, we divide the sum over $r$ as $0<r \leq \xi$ and $\xi < r \leq Ax^{1/2}$ for some $\xi > 0$ chosen explicitly later. Then the first sum over $0 < r \leq \xi$ is evaluated as 
\begin{align}\label{first-xi}
	\left\lvert \sum_{0 < r \leq \xi} \sum_{s} e^{-2\pi i m\overline{s}/r} \varphi_m(r,s) \right\rvert \leq \sum_{0 < s \leq r \leq \xi} 1 = O(\xi^2).
\end{align}
As for the second part, we proceed as if the set $I$ is an interval $[L,R]$. We can evaluate each interval similarly if it is a union of two intervals. For the inner sum, by partial summation formula, we have
\begin{align*}
	\sum_{\substack{L \leq s \leq R \\ (r,s) = 1}} e^{-2\pi i m\overline{s}/r} \varphi_m(r,s) = g(R) \varphi_m(r, R) - \sum_{\mu=L}^{R-1} g(\mu) (\varphi_m(r,\mu+1) - \varphi_m(r, \mu)),
\end{align*}
where we put
\[
	g(\mu) = \sum_{\substack{L \leq s \leq \mu \\ (r,s) = 1}} e^{-2\pi i m\overline{s}/r}.
\]
Since $\mu - L \leq r$, we can apply the estimate on an incomplete Kloosterman sum given by Hooley~\cite[Lemma 3]{Hoo63}, and we have
\[
	g(\mu) = O\left(r^{1/2} (\gcd(m,r))^{1/2} d(r) \log(2r) \right),
\]
where $d(r)$ is the number of positive divisors of $r$. Moreover, by the mean value theorem and the conditions $0 < s \leq r$ and $cr^2 + (d-a)rs - bs^2 \geq Br^2$, we see that
\[
	\varphi_m(r,\mu+1) - \varphi_m(r,\mu) = O \left(\frac{m}{r^3} \right).
\]
Therefore, we get
\begin{align*}
	&\sum_{\substack{L \leq s \leq R \\ (r,s) = 1}} e^{-2\pi i m\overline{s}/r} \varphi_m(r,s)\\
	&= O\left(r^{1/2} (\gcd(m,r))^{1/2} d(r) \log(2r) \right) + O\left(m r^{-3/2} (\gcd(m,r))^{1/2} d(r) \log(2r) \right).
\end{align*}

Next, we consider the sum over $r$. For the first term, we have
\begin{align*}
	\sum_{\xi < r \leq Ax^{1/2}} r^{1/2} (\gcd(m,r))^{1/2} d(r) \log(2r) &\ll x^{1/4} \log x \sum_{1 \leq r \leq Ax^{1/2}} (\gcd(m,r))^{1/2} d(r)\\
		&= O \left(\sigma_{-1/2}^2(m) x^{3/4} \log^2 x \right),
\end{align*}
where we used~\cite[Lemma 4]{Hoo63}, namely
\begin{align}\label{Hooley-Lemma-4}
	\sum_{r \leq x} (\gcd(m,r))^{1/2} d(r) = O \left(x \log(2x) \sigma_{-1/2}^2(m) \right).
\end{align}
As for the second term, by partial summation formula and \eqref{Hooley-Lemma-4} again, 
\begin{align*}
	&\sum_{\xi < r \leq Ax^{1/2}} m r^{-3/2} (\gcd(m,r))^{1/2} d(r) \log(2r)\\
	&\ll m \log x \sum_{\xi < r \leq A x^{1/2}} (\gcd(m,r))^{1/2} d(r) \cdot r^{-3/2} \ll m \sigma_{-1/2}^2(m) \xi^{-1/2} \log^2 x.
\end{align*}
	
Therefore, the second part of \eqref{Sgamma-expre} is evaluated as
\begin{align*}
	&\sum_{\xi < r \leq A x^{1/2}} \sum_{\substack{L \leq s \leq R \\ (r,s) = 1}} e^{-2\pi i m\overline{s}/r} \varphi_m(r,s) = O_\epsilon \left(m^\epsilon x^{3/4+\epsilon} \right) + O_\epsilon \left(m^{1+\epsilon} \xi^{-1/2} x^\epsilon \right).
\end{align*}
Combining it with the evaluation~\eqref{first-xi} of the first part implies
\[
	\sum_{0 < k \leq x} S_\gamma(k,m) = O(\xi^2) + O_\epsilon \left(m^\epsilon x^{3/4+\epsilon} \right) + O_\epsilon \left(m^{1+\epsilon} \xi^{-1/2} x^\epsilon \right).
\]
If we choose $\xi = m^{2/5}$, it holds that $\xi^2 = m \xi^{-1/2} = m^{4/5}$, which completes the proof of \cref{Hooley-theorem1}.

\section{Eisenstein series of weight $0$} \label{Appendix-C}

In this last section, we note one remarkable result on the limit formula for the weight $0$ analogue of the Eisenstein series. Similar as in \cref{s3}, for each non-scalar element $\gamma \in \Gamma$, we define the real analytic Eisenstein series by
\[
	E_{0,\gamma}(z,s) := \sum_{Q \sim Q_\gamma} \frac{y^s}{\lvert Q(z,1)\rvert^s} = \sum_{g \in \Gamma_{w_\gamma} \backslash \Gamma} \frac{y^s}{\lvert Q_{g^{-1} \gamma g}(z,1)\rvert^s}, \qquad \Re(s) > 1,
\]
without the sign function. This series converges absolutely and locally uniformly for $\Re(s) > 1$ and $z \in \bbH$, and satisfies $E_{0, \gamma}(gz, s) = E_{0, \gamma}(z,s)$ for any $g \in \Gamma$. As in \cref{s3-1}, for $\gamma = T = \smat{1& 1 \\ 0 & 1}$, this is the classical Eisenstein series
\[
	E_{0, T}(z,s) = \sum_{g \in \Gamma_\infty \backslash \Gamma} \Im(gz)^s = \frac{1}{2} \sum_{(c,d) = 1} \frac{y^s}{\lvert cz+d\rvert^{2s}}.
\]
On the other hand, for elliptic or hyperbolic elements, we have
\begin{align}\label{E-ell-hyp-0}
	\begin{split}
	E_{0,\gamma}(z, s) = \left\{\begin{array}{ll}
		 \displaystyle{\lvert D_\gamma\rvert^{-s/2} \sum_{g \in \Gamma_{w_\gamma} \backslash \Gamma} \frac{1}{\sinh(\hyp{g z, w_\gamma})^s}} \quad &\text{if } \gamma \text{ is elliptic},\\
		\displaystyle{D_\gamma^{-s/2} \sum_{g \in \Gamma_{w_\gamma} \backslash \Gamma} \frac{1}{\cosh(\hyp{g z, S_\gamma})^s}} &\text{if } \gamma \text{ is hyperbolic},
	\end{array} \right.
	\end{split}
\end{align}
where we denote by $d_{\mathrm{hyp}}$ the metric on the upper half-plane $\bbH$. This metric is given by
\begin{align}\label{dis-def}
	\hyp{z, \tau} = \arccosh \left(1 + \frac{\lvert z - \tau\rvert^2}{2\Im(z) \Im(\tau)} \right) = \arcsinh \left(\frac{\lvert (z-\tau)(z- \bar{\tau})\rvert}{2\Im(z) \Im(\tau)} \right),
\end{align}
and satisfies $\hyp{gz, g\tau} = \hyp{z, \tau}$ for any $g \in \SL_2(\mathbb{R})$. The expression \eqref{E-ell-hyp-0} immediately follows from the next lemma.

\begin{lem}\label{dis-hyp}
	Let $\gamma \in \Gamma$ be a hyperbolic element, and $D_\gamma := \tr(\gamma)^2 - 4$ the discriminant of $Q_\gamma(X, Y)$. Then we have
	\[
		\lvert Q_\gamma(z,1)\rvert = \sqrt{D_\gamma} \Im(z) \cosh(\hyp{z, S_\gamma}),
	\]
	where $S_\gamma$ is the geodesic in $\bbH$ connecting two fixed points $w_\gamma > w'_\gamma$ of $\gamma$.
\end{lem}

\begin{proof}
	This lemma is described in~\cite[Lemma 2.5.4]{Vol18}. We now review his proof. Let $\calI$ be the $y$-axis. By~\cite[Theorem 7.9.1 (ii)]{Bea83}, we have
	\begin{align*}
		\cosh(\hyp{z, S_\gamma}) &= \cosh(\hyp{M_\gamma^{-1}z, \calI}) = \frac{1}{\sin \arg(M_\gamma^{-1} z)} = \frac{\lvert M_\gamma^{-1} z \rvert}{\Im(M_\gamma^{-1} z)}\\
			&= \frac{\lvert Q_{M_\gamma^{-1} \gamma M_\gamma}(M_\gamma^{-1} z ,1)\rvert}{\sqrt{D_\gamma} \Im(M_\gamma^{-1} z)} = \frac{1}{\sqrt{D_\gamma}} \frac{\lvert Q_\gamma(z,1)\rvert}{\Im(z)}.
	\end{align*}
\end{proof}

The hyperbolic case was originally studied by Petersson~\cite{Pet44} and Kudla--Millson~\cite{KM79} for weight $2$. More recently, Jorgenson--Kramer--von Pippich, and so on~\cite{JKvP10,JvPS16,Pippich2016,vPSV17,Vol18} studied the above weight $0$ analogue.

Our goal in this section is getting the limit formulas for these Eisenstein series. In the classical case, the Eisenstein series $E_{0, T}(z,s)$ is meromorphically continued to the whole $s$-plane, and we know that
\begin{align}\label{par-Eisen}
	\begin{split}
		E_{0,T}(z,s) &= 1 + \log(y \lvert \eta(z)\rvert^4) \cdot s + O(s^2)\\
			&= \frac{3/\pi}{s-1} - \frac{3}{\pi} \log(y \lvert \eta(z) \rvert^4) + \frac{6}{\pi} \bigg( \gamma - \log 2 - \frac{6 \zeta'(2)}{\pi^2} \bigg) + O(s-1).
	\end{split} 
\end{align}
This is called the Kronecker limit formula, where $\eta(z) = q^{1/24} \prod_{n=1}^\infty (1-q^n)$ is the Dedekind eta function, $\gamma$ is the Euler constant, and $\zeta(s)$ is the Riemann zeta function. (We hope there is no confusion between the notation $\gamma$ used for a matrix and the Euler constant). In the elliptic case, von Pippich~\cite{Pippich2016} gave the meromorphic continuation of the elliptic Eisenstein series $E_{0,\gamma} (z,s)$ for $\gamma = S = \smat{0 & -1 \\ 1 & 0}, U = \smat{1 & -1 \\ 1 & 0}$, and established the limit formula
\[
	E_{0,\gamma} (z,s) = - \frac{2}{\lvert \Gamma_{w_\gamma}\rvert} \log \lvert j(z) - j(w_\gamma) \rvert \cdot s + O(s^2).
\]
Here the number $\lvert \Gamma_{w_\gamma}\rvert$ is given by $4$ or $6$ according as $\gamma = S$ or $U$. On the other hand, for a hyperbolic $\gamma$, Jorgenson--Kramer--von Pippich~\cite{JKvP10} derived the meromorphic continuation of the hyperbolic Eisenstein series $E_{0,\gamma} (z,s)$ by using the spectral expansion. According to their results, the function $E_{0,\gamma} (z,s)$ has a double zero at $s=0$. Moreover, von Pippich--Schwagenscheidt--V\"{o}lz~\cite[Remark 5.7 in arXiv version]{vPSV17} described that it is an interesting problem to investigate the second order coefficient of $E_{0, \gamma} (z,s)$ at $s=0$.

Our following result provides an answer to this problem. Similar as before, we can easily reduce the problem to the case that $\gamma$ is a primitive hyperbolic element with $\sgn(Q_\gamma) > 0$ and $\tr(\gamma) > 2$.

\begin{thm}\label{hyp-limit}
	Let $\gamma = \smat{a & b \\ c & d} \in \Gamma$ be a primitive hyperbolic element with $\sgn(Q_\gamma) = \sgn(c) > 0$ and $\tr(\gamma) = a+d > 2$. Then the function $E_{0, \gamma}(z,s)$ is analytically continued to $s=0$, and equals
	\begin{align*}
		&E_{0, \gamma} (z,s)\\
			& = -\frac{1}{2} \left( \int_{\tau_0}^{\gamma \tau_0} \log \lvert j(z) - j(\tau)\rvert \frac{-\sqrt{D_\gamma} d\tau}{Q_\gamma(\tau, 1)} + \sum_{g \in \Gamma_{w_\gamma} \backslash \Gamma} \arcsin^2 \left(\frac{1}{\cosh(\hyp{gz, S_\gamma})} \right) \right) s^2 + O(s^3),
	\end{align*}
	where $\tau_0 \in S_\gamma$ is any point on $S_\gamma$, and the path of integration is on the geodesic $S_\gamma$.
\end{thm}

\begin{rem}
	It is also an interesting problem to investigate the Poincar\'{e} series appeared in the above theorem. It is known that the hyperbolic Eisenstein series $E_{0,\gamma}(z,s)$ satisfies the differential equation
	\[
		(\Delta_0 -s(1-s)) E_{0, \gamma}(z,s) = s^2 D_\gamma E_{0, \gamma}(z, s+2).
	\]
	Here $\Delta_0$ is the hyperbolic Laplacian
	\[
		\Delta_0 := -y^2 \left(\frac{\partial^2}{\partial x^2} + \frac{\partial^2}{\partial y^2} \right) = -\sin^2 \theta \left(\frac{\partial^2}{\partial \rho^2} + \frac{\partial^2}{\partial \theta^2} \right)
	\]
	for $z = x+iy = e^{\rho + i\theta}$. We put $E_{0,\gamma}(z,s) = \calL_\gamma(z) s^2 + O(s^3)$. By this differential equation, we have
	\begin{align}\label{Delta-image}
		\Delta_0 \calL_\gamma(z) = D_\gamma E_{0,\gamma}(z,2) = \sum_{g \in \Gamma_{w_\gamma} \backslash \Gamma} \frac{1}{\cosh(d_{\mathrm{hyp}}(gz, S_\gamma))^2} = \sum_{g \in \Gamma_{w_\gamma} \backslash \Gamma} (\sin \arg(M_\gamma^{-1} z))^2.
	\end{align}
	Therefore, the Poincar\'{e} series in \cref{hyp-limit} is a natural preimage of \eqref{Delta-image} under $\Delta_0$. The first term in the theorem is the harmonic part, which is annihilated by $\Delta_0$. 
\end{rem}

\begin{rem}
	Lagarias--Rhoades~\cite{LR16} and the author~\cite{Mat19,Mat20} studied higher-order Laurent coefficients of the classical parabolic Eisenstein series $E_{0,T}(z,s)$ at $s = 0$, which are called polyharmonic Maass forms. More recently, Bringmann--Kane~\cite{BK20} introduced polar polyharmonic Maass forms in terms of the Laurent coefficients of the automorphic Green function $G(z,\tau; s)$ defined below. One natural question is whether there exists a suitable theory containing the above limit formula for $E_{0, \gamma}(z,s)$.
\end{rem}

Here is what we do. To prove this theorem, we first generalize von Pippich's result. For any $\Gamma$-inequivalent $z, \tau \in \bbH$, we consider the function
\[
	P(z,\tau; s) := \sum_{g \in \Gamma} \frac{1}{\sinh(\hyp{gz, \tau})^s}, \qquad \Re(s) > 1,
\]
and show that
\begin{align}\label{Ew-Limit}
	P(z,\tau; s) = -2 \log \lvert j(z) - j(\tau)\rvert \cdot s + O(s^2)
\end{align}
after the meromorphic continuation to $s = 0$. Actually this limit formula follows from the results of Gross--Zagier~\cite{GZ85} and Fay~\cite{Fay77}. The proofs of \eqref{Ew-Limit} and \cref{hyp-limit} are described by means of the intermediate function
\[
	Q(z,\tau; s) := \sum_{g \in \Gamma} \frac{1}{\cosh(\hyp{gz,\tau})^s}, \qquad \Re(s) > 1.
\]

To make this precise, we now recall the results of Gross-Zagier and Fay. For $\Re(s) > 1$ and any $\Gamma$-inequivalent $z, \tau \in \bbH$, we consider the automorphic Green function defined by
\[
	G(z,\tau; s) := -\frac{2^s \Gamma(s)^2}{\Gamma(2s)} \sum_{g \in \Gamma} \bigg(1 + \cosh(\hyp{gz, \tau}) \bigg)^{-s} {}_2 F_1 \bigg(s, s; 2s; \frac{2}{1 + \cosh(\hyp{gz,\tau})} \bigg),
\]
where ${}_2 F_1(a, b; c; z)$ is Gauss' hypergeometric function given by
\[
	{}_2 F_1(a,b; c; z) := \sum_{k=0}^\infty \frac{(a)_k (b)_k}{(c)_k} \frac{z^k}{k!}
\]
with the Pochhammer symbol $(x)_k := x(x+1) \cdots (x+k-1)$. We note that Gross--Zagier's definition given in~\cite{GZ85} equals $G(z,\tau;s)/2$, and Fay's definition in~\cite{Fay77} equals $G(z,\tau;s)/8\pi$. 

\begin{prop}~\cite[Proposition\ 5.1]{GZ85}
	For above $z, \tau \in \bbH$, the function $G(z,\tau;s)$ has the meromorphic continuation to the whole $s$-space, and we have
	\begin{align*}
		&4\log \lvert j(z) - j(\tau)\rvert + 48\\
			&= \lim_{s \to 1} \bigg( G(z,\tau;s) + 8\pi E_{0,T}(z,s) + 8\pi E_{0,T}(\tau,s) - 8\pi \frac{\Gamma \left(\frac{1}{2} \right) \Gamma \left(s - \frac{1}{2} \right)}{\Gamma(s)} \frac{\zeta(2s-1)}{\zeta(2s)} \bigg).
	\end{align*}
\end{prop}

By combining this proposition and the Kronecker limit formula \eqref{par-Eisen}, we have
\[
	G(z,\tau;s) = \frac{-24}{s-1} + f(z,\tau) + O(s-1),
\]
where $\tau = u + iv$ and 
\[
	f(z,\tau) = 4\log \lvert j(z) - j(\tau)\rvert + 24 \log (y \lvert \eta(z) \rvert^4) + 24 \log(v \lvert \eta(\tau)\rvert^4) -48 \bigg(\gamma - \log 2 - \frac{6\zeta'(2)}{\pi^2} -1 \bigg).
\]

\begin{prop}~\cite[(79)]{Fay77}
	For above $z, \tau \in \bbH$, the function $G(z,\tau;s)$ satisfies the functional equation
	\[
		G(z,\tau;s) - G(z,\tau;1-s) = \frac{8\pi}{1-2s} E_{0,T}(z, 1-s) E_{0,T}(\tau,s).
	\]
\end{prop}

By this functional equation, we obtain the Taylor expansion at $s = 0$,
\begin{align}\label{G-Tay}
	G(z,\tau;s) = 4 \log \lvert j(z) - j(\tau)\rvert + O(s).
\end{align}
To show \eqref{Ew-Limit}, we consider the intermediate function $Q(z,\tau;s)$. On this function, we have the following two lemmas, which are described in von Pippich's preprint~\cite{Pippich2016}.

\begin{lem}\label{E=P}
	For $\Gamma$-inequivalent $z,\tau \in \bbH$ and $\Re(s) > 1$, we have
	\[
		P(z,\tau; s) = \sum_{k=0}^\infty \frac{\left(\frac{s}{2}\right)_k}{k!} Q(z,\tau; s+2k),
	\]
	which converges absolutely and locally uniformly for $\Re(s) > 1$ and $z \in \bbH \setminus \Gamma \tau$.
\end{lem}

\begin{proof}
	We first check the convergence of the right-hand side. For fixed $\tau \in \bbH$ and $s$ with $\sigma := \Re(s) > 1$, we have
	\begin{align*}
		\left| \sum_{k=0}^\infty \frac{\left(\frac{s}{2}\right)_k}{k!} Q(z,\tau; s+2k) \right| &\leq \sum_{k=0}^\infty \frac{\lvert \left(\frac{s}{2} \right)_k \rvert}{k!} \sum_{g \in \Gamma} \frac{1}{\cosh(d_{\mathrm{hyp}}(gz, \tau))^{\sigma+2k}}.
	\end{align*}
	By the trivial bound $\lvert \Gamma(s) \rvert \leq \Gamma(\sigma)$, we get $\lvert (s/2)_k\rvert = \lvert \Gamma(s/2+k)/\Gamma(s/2) \rvert \leq \Gamma(\sigma/2+k)/\lvert \Gamma(s/2)\rvert$. Also we get the bound $\cosh(d_{\mathrm{hyp}}(gz, \tau)) \geq C > 1$ for any compact subset $z \in K \subset \bbH \setminus \Gamma \tau$, where the constant $C$ depends only on $K$ and $\tau$. Thus we obtain
	\begin{align*}
		\left| \sum_{k=0}^\infty \frac{\left(\frac{s}{2}\right)_k}{k!} Q(z,\tau; s+2k) \right| &\leq \frac{\Gamma \left(\frac{\sigma}{2} \right)}{\lvert\Gamma \left( \frac{s}{2} \right) \rvert} \sum_{g \in \Gamma} \frac{1}{\cosh(d_{\mathrm{hyp}}(gz, \tau))^\sigma} \sum_{k=0}^\infty \frac{\left(\frac{\sigma}{2} \right)_k}{k!} \frac{1}{C^{2k}}\\
			&= \frac{\Gamma \left(\frac{\sigma}{2} \right)}{\lvert\Gamma \left( \frac{s}{2} \right) \rvert} \left(1 - \frac{1}{C^2} \right)^{-\frac{\sigma}{2}} Q(z,\tau; \sigma) < \infty.
	\end{align*}
	Similarly for a fixed $z \in \bbH \setminus \Gamma \tau$, we can show the locally uniformly convergence of this series for $\Re(s) > 1$. Then we can change the order of sums in the right-hand side. Finally the equation
	\[
		\sum_{k=0}^\infty \frac{\left(\frac{s}{2}\right)_k}{k!} \frac{1}{\cosh(\hyp{gz, \tau})^{s + 2k}} = \frac{1}{\sinh(\hyp{gz, \tau})^s}
	\]
	concludes the proof.
\end{proof}

\begin{lem}
	For $\Gamma$-inequivalent $z, \tau \in \bbH$ and $\Re(s) > 1$, we have
	\[
		G(z,\tau;s) = - \frac{2^s \Gamma(s)^2}{\Gamma(2s)} \sum_{k=0}^\infty \frac{\left( \frac{s}{2} \right)_k \left(\frac{s+1}{2} \right)_k}{ \left(s + \frac{1}{2} \right)_k k!} Q(z,\tau; s+2k).
	\]
\end{lem}

\begin{proof}
	This follows from the formula given in~\cite[9.134.1]{GR07},
	\[
		{}_2 F_1 \left(\alpha, \beta; 2\beta; z \right) = \left(1 - \frac{z}{2} \right)^{-\alpha} {}_2 F_1 \left(\frac{\alpha}{2}, \frac{\alpha+1}{2}; \beta + \frac{1}{2}; \left(\frac{z}{2-z} \right)^2 \right).
	\]
	The details are similar as in the proof of \cref{E=P}.
\end{proof}

By using the lemmas, we easily see that the right-hand side of the function
\[
	P(z,\tau; s) + \frac{\Gamma(2s)}{2^s \Gamma(s)^2} G(z,\tau;s) = \sum_{k=1}^\infty \frac{\left(\frac{s}{2} \right)_k}{k!}  \left(1- \frac{\left(\frac{s+1}{2} \right)_k}{\left(s+\frac{1}{2} \right)_k} \right) Q(z, \tau; s+2k)
\]
converges in $\Re(s) > -1/2$, and has a double zero at $s=0$. From \eqref{G-Tay}, we obtain \eqref{Ew-Limit}. \\

Next we consider the relation between the intermediate function $Q(z,\tau; s)$ and the hyperbolic Eisenstein series $E_{0,\gamma}(z,s)$. The key object is the cycle integral
\[
	\int_{\tau_0}^{\gamma \tau_0} Q(z, \tau; s) \frac{d\tau}{Q_\gamma(\tau, 1)}.
\]
Here the path of integration is in the geodesic $S_\gamma$. Since $Q(z, \tau; s)$ is $\Gamma$-invariant in $\tau$, this integral is independent of the choice of $\tau_0 \in S_\gamma$. By the unfolding argument, for $\Re(s) > 1$,
\begin{align*}
	\int_{\tau_0}^{\gamma \tau_0} Q(z,\tau; s) \frac{d\tau}{Q_\gamma(\tau,1)} &= 2 \int_{\tau_0}^{\gamma \tau_0} \sum_{g \in \Gamma_{w_\gamma} \backslash \Gamma} \sum_{n \in \mathbb{Z}} \frac{1}{\cosh(\hyp{\gamma^n gz, \tau})^s} \frac{d\tau}{Q_\gamma(\tau, 1)}\\
		&= 2 \sum_{g \in \Gamma_{w_\gamma} \backslash \Gamma} \int_{S_\gamma} \frac{1}{\cosh(\hyp{gz, \tau})^s} \frac{d\tau}{Q_\gamma(\tau, 1)},
\end{align*}
where $S_\gamma$ is oriented from $w'_\gamma$ to $w_\gamma$. By~\cite[Theorem 7.11.1]{Bea83}, we have
\[
	\cosh(\hyp{gz, \tau}) = \cosh(\hyp{M_\gamma \cdot i \lvert M_\gamma^{-1} gz \rvert, \tau}) \cosh(\hyp{gz, S_\gamma}).
\]
For any $\z \in \bbH$, we see that
\begin{align*}
	\int_{S_\gamma} \frac{1}{\cosh(\hyp{M_\gamma i \lvert \z \rvert, \tau})^s} \frac{d\tau}{Q_\gamma(\tau, 1)} &= -\frac{1}{\sqrt{D_\gamma}} \int_{-\infty}^\infty \frac{dr}{\cosh(\hyp{i \lvert \z\rvert, ie^r})^s}\\
		&= -\frac{1}{\sqrt{D_\gamma}} \int_{-\infty}^\infty \frac{dr}{\cosh(r - \log \lvert \z \rvert)^s}\\
		&= -\frac{1}{\sqrt{D_\gamma}} \int_{-\infty}^\infty \frac{dr}{\cosh(r)^s}.
\end{align*}
From the formula~\cite[p.10]{MOS66}
\[
	\int_0^\infty \sinh(t)^\alpha \cosh(t)^{-\beta} dt = \frac{\Gamma \left(\frac{\alpha +1}{2} \right) \Gamma \left(\frac{\beta- \alpha}{2} \right)}{2 \Gamma \left(\frac{\beta +1}{2} \right)}, \qquad (\Re(\beta - \alpha) > 0, \Re(\alpha) > -1),
\]
this integral equals
\[
	\int_{S_\gamma} \frac{1}{\cosh(\hyp{M_\gamma i \lvert \z \rvert, \tau})^s} \frac{d\tau}{Q_\gamma(\tau, 1)} = -\frac{1}{\sqrt{D_\gamma}} \frac{\sqrt{\pi} \Gamma \left(\frac{s}{2} \right)}{\Gamma \left(\frac{s+1}{2} \right)}.
\]
Therefore for $\Re(s) > 1$, we obtain that
\begin{align*}
	\int_{\tau_0}^{\gamma \tau_0} Q(z,\tau; s) \frac{d\tau}{Q_\gamma(\tau,1)} &= -\frac{2}{\sqrt{D_\gamma}} \frac{\sqrt{\pi} \Gamma \left(\frac{s}{2} \right)}{\Gamma \left(\frac{s+1}{2} \right)} \sum_{g \in \Gamma_{w_\gamma} \backslash \Gamma} \frac{1}{\cosh(\hyp{gz, S_\gamma})^s} \\
		&= -\frac{2D_\gamma^{\frac{s-1}{2}} \sqrt{\pi} \Gamma \left(\frac{s}{2} \right)}{\Gamma \left(\frac{s+1}{2} \right)} E_{0, \gamma}(z,s).
\end{align*}

\begin{thm}\label{E=intP}
	Let $\gamma$ be a primitive hyperbolic element with $\sgn(Q_\gamma) > 0$ and $\tr(\gamma) > 2$. For $\Re(s) > 1$, we have
	\[
		E_{0, \gamma}(z,s) = -\frac{\Gamma \left(\frac{s+1}{2} \right)}{2 D_\gamma^{\frac{s-1}{2}} \sqrt{\pi} \Gamma \left(\frac{s}{2} \right)} \int_{\tau_0}^{\gamma \tau_0} Q(z, \tau; s) \frac{d\tau}{Q_\gamma(\tau, 1)}.
	\]
\end{thm}

By using these results, we now prove \cref{hyp-limit}. For $\Re(s) > 1$, we have
\[
	Q(z,\tau; s) = P(z,\tau; s) - \sum_{k=1}^\infty \frac{\left(\frac{s}{2} \right)_k}{k!} Q(z,\tau; s+2k)
\]
by \cref{E=P}. By plugging this in \cref{E=intP},
\[
	E_{0, \gamma} (z,s) =  -\frac{\Gamma \left(\frac{s+1}{2} \right)}{2 D_\gamma^{\frac{s-1}{2}} \sqrt{\pi} \Gamma \left(\frac{s}{2} \right)} \int_{\tau_0}^{\gamma \tau_0} \left( P(z,\tau; s) - \sum_{k=1}^\infty \frac{\left(\frac{s}{2} \right)_k}{k!} Q(z,\tau; s+2k) \right) \frac{d\tau}{Q_\gamma(\tau, 1)}
\]
Since the sum in the integral converges absolutely and locally uniformly, the termwise integration is legitimate when $z$ is not on the net of geodesics $S_{g^{-1} \gamma g}$. Then the expression
\begin{align*}
	E_{0, \gamma} (z,s) &=  -\frac{\Gamma \left(\frac{s+1}{2} \right)}{2 D_\gamma^{\frac{s-1}{2}} \sqrt{\pi} \Gamma \left(\frac{s}{2} \right)} \left( \int_{\tau_0}^{\gamma \tau_0} P(z,\tau; s) \frac{d\tau}{Q_\gamma(\tau,1)} - \sum_{k=1}^\infty \frac{\left(\frac{s}{2} \right)_k}{k!} \int_{\tau_0}^{\gamma \tau_0} Q(z,\tau; s+2k) \frac{d\tau}{Q_\gamma(\tau,1)} \right) \\
		&= -\frac{\Gamma \left(\frac{s+1}{2} \right)}{2 D_\gamma^{\frac{s-1}{2}} \sqrt{\pi} \Gamma \left(\frac{s}{2} \right)} \int_{\tau_0}^{\gamma \tau_0} P(z,\tau; s) \frac{d\tau}{Q_\gamma(\tau,1)} - \sum_{k=1}^\infty \frac{D_\gamma^k \left(\frac{s}{2} \right)_k^2}{\left(\frac{s+1}{2} \right)_k k!} E_{0, \gamma}(z, s+2k)
\end{align*}
gives the meromorphic continuation of $E_\gamma(z,s)$ to $s = 0$. Actually, the infinite sum in the right-hand side converges absolutely and locally uniformly for $\Re(s) > -1$. At $s=0$ we have the Taylor expansion of the form
\begin{align*}
	E_{0,\gamma}(z,s) = \left( \frac{\sqrt{D_\gamma}}{2} \int_{\tau_0}^{\gamma \tau_0} \log \lvert j(z) - j(\tau) \rvert \frac{d\tau}{Q_\gamma(\tau,1)} - \frac{1}{4} \sum_{k=1}^\infty \frac{\sqrt{\pi} D_\gamma^k \Gamma(k)}{k \Gamma \left(k + \frac{1}{2} \right)} E_{0, \gamma}(z,2k) \right) s^2 + O(s^3).
\end{align*}
For the second term, we easily see that
\begin{align*}
	\frac{1}{4} \sum_{k=1}^\infty \frac{\sqrt{\pi} D_\gamma^k \Gamma(k)}{k \Gamma \left(k+\frac{1}{2} \right)} E_{0,\gamma} (z, 2k) &= \frac{1}{4} \sum_{k=1}^\infty \frac{\sqrt{\pi} \Gamma(k)}{k \Gamma \left(k+\frac{1}{2} \right)} \sum_{g \in \Gamma_{w_\gamma} \backslash \Gamma} \frac{1}{\cosh(\hyp{gz, S_\gamma})^{2k}}\\
		&= \frac{1}{2} \sum_{g \in \Gamma_{w_\gamma} \backslash \Gamma} \arcsin^2 \left(\frac{1}{\cosh(d_{\mathrm{hyp}}(gz, S_\gamma))} \right).
\end{align*}
Therefore we finally obtain \cref{hyp-limit}.

\section*{Acknowledgements} 
It is a pleasure to thank \"{O}zlem Imamo\={g}lu for introducing the author to the hyperbolic Rademacher symbol, and Markus Schwagenscheidt for some enlightening discussions. The author also would like to express his unbounded gratitude to Masanobu Kaneko for offering his wisdom and encouragement. This work was supported by JSPS KAKENHI Grant Numbers JP18J20590, JP20K14292, and JP21K18141.


\bibliographystyle{amsplain}
\bibliography{references}


\end{document}